\documentclass[11pt]{article}
\usepackage[english,activeacute]{babel}
\usepackage[utf8]{inputenc}
\usepackage[Lenny]{fncychap}
\usepackage[hidelinks]{hyperref}
\usepackage{amssymb,amsmath,amsthm} 
\usepackage{graphics,graphicx} 
\usepackage{subfigure} 
\DeclareGraphicsExtensions{.pdf,.png,.jpg,.eps}
\usepackage{lscape}
\usepackage{amsfonts,mathrsfs}
\usepackage{multicol}
\usepackage{fancybox,calc}

\usepackage{appendix}

\usepackage[usenames,dvipsnames]{color}

\usepackage{makeidx} 
\usepackage{appendix} 

\newtheorem{Obs}{Remark}[section]

\newtheorem{teo}{Theorem}[section]
\newtheorem{proposition}{Proposition}[section]

\newtheorem{lema}{Lemma}[section]
\newtheorem{Definition}{Definition}[section]

\numberwithin{figure}{section}
\numberwithin{table}{subsection}
\numberwithin{table}{section}
\numberwithin{equation}{section}

\numberwithin{equation}{section}

\makeindex 
\begin{document}
\title{Robust Stackelberg controllability for the Navier--Stokes equations}
\author{Cristhian Montoya and Luz de Teresa\thanks{Universidad Nacional Aut\'onoma de M\'exico, UNAM, Instituto de 
    Matem\'aticas,  Ciudad de M\'exico, M\'exico.
    emails: cmontoya@matem.unam.mx, \quad ldeteresa@im.unam.mx. The research was partially supported by UNAM-DGAPA-PAPIIT
IN102116. The first author is supported by project FORDECYT 265667 of CONACYT}}
\date{}
\maketitle
\abstract{In this paper we deal with a robust Stackelberg strategy for the  
    Navier--Stokes system. The scheme is based in considering a robust control problem
    for the ``follower control'' and its associated disturbance function. Afterwards, we consider the 
    notion of Stackelberg optimization (which is associated to the ``leader control'') 
    in order to deduce a local null controllability result for the Navier--Stokes system.}
MSC 35Q30, 35Q93, 49J20, 91A65, 93C10  


\section{\normalsize Introduction}
	The theory of robust control began in the late 1970s and early 
    1980s for finite dimensional systems. Since then,  many techniques have been developed to deal with
    systems with uncertainties. In the late 90s the papers of Bewley et. al 
    \cite{bewley2000existence} presented 
    the first rigorous generalization of the concepts in the case of partial differential equations. 
    What could we understand by  robustness in a control system?
    Well, informally, a 
    controller designed for a particular set of parameters is said to be robust if it also functions 
 	correctly under a uncertainty: the controller is designed to work assuming 
    that certain variable will be unknown. In this sense, one could think in the 
    worst--case disturbance of the  system, and design a controller which is suited to handle even this 
    extreme situation. Thus, the problem of finding a robust control involves the problem of finding
    the worst--case disturbance in the spirit of a non--cooperative game (when there is not 
    cooperation between the controller and disturbance function), which is from the 
    mathematical point of view to reach a saddle point for the pair disturbance--controller.

    The research on robust control for PDE systems  is in an early stage. Much of the literature deals
    with numerical aspects and much of the theory has been developed for fluid mechanics and for some 
    elliptic problems. See e.g. \cite{abergel1990some, bewley2000existence, bewley1997optimal,
    medjo2008optimal, herzog2012weak}. In this paper we will present 
    a hierarchic strategy to deal  with robust control and, simultaneously, with null control for 
   	incompressible fluids  modelled by the Navier--Stokes equations with Dirichlet boundary conditions. 
   	
    We will work in the setting of  a Stackelberg competition, see \cite{von1952theory}. This consists in a non--
	cooperative decision problem in which one of the participants enforce its strategy on the other participants. 
    We assume that we can act on the dynamics of the system through a hierarchy of controls. In our case
    the controls are external forces acting on the system, where the leader control has a local null 
    controllability objective while the follower control and perturbation solve a robust control problem.

 	To be precise: let $\Omega$  be a nonempty bounded connected open subset of 
    $\mathbb{R}^N$ ($N=2$ or $N=3$) of class $C^{\infty}$. Let $T > 0$ and 
    let   $\omega$ and $\mathcal{O}$ be (small) nonempty open subsets of $\Omega$ 
    with $\omega\cap \mathcal{O}=\emptyset$. We will use the notation 
    $Q:=\Omega\times(0,T)$,\, $\Sigma:=\partial\Omega\times(0,T)$ and   $n(x)$   will denote  
    the outward unit normal vector  at the point $x\in\partial\Omega$.\\
    Let us consider the Navier-Stokes system with homogeneous Dirichlet boundary
    conditions
\begin{equation}\label{1.1.main_system}
    \left\{
    \begin{aligned}
    \begin{array}{llll}
    y_{t}-\Delta y+(y\cdot\nabla)y+\nabla p=h1_{\omega}+v\chi_{\mathcal{O}}+\psi  &
    \text{ in }& Q,\\
    \nabla \cdot y=0&  \text{ in }& Q, \\
    y=0&\text{ on }&\Sigma, \\
    y(\cdot,0)=y_0(\cdot)& \text{ in }&\Omega,
    \end{array}
    \end{aligned}\right.
\end{equation}
    where $h=h(x,t)\in L^2(\omega\times(0,T))^N$ is called the ``leader control'', 
    $v=v(x,t)\in L^2(\mathcal{O}\times(0,T))^N$  is the ``follower control'', 
    $\psi\in L^2(Q)^N$ is an unknown perturbation  and $y_0$ an initial state in a suitable space.
    Here $1_{\omega}$ is the characteristic function of the set $\omega$ and $\chi_{\mathcal{O}}$ is a smooth
	non-negative function such that $supp\,\,\chi_\mathcal{O}=\overline{\mathcal{O}}$.

    To our knowledge there are not results in the literature concerning a robust Stackelberg strategy for 
    system \eqref{1.1.main_system}. As far as we know, the first paper  on robust
    Stackelberg controllability is \cite{hernandez2016robust}, which develops the concept of control for a 
    semi--linear parabolic equation.  However, there exist several 
	papers which treat independently robust and hierarchical control for the 
    Navier--Stokes system.  In the context of robust control 
    (that is  $h\equiv 0$ in \eqref{1.1.main_system}), 
    the works \cite{bewley2000general} and \cite{bewley2000existence} 
    show  the existence and uniqueness of the solution to the robust control problem for the $N$-dimensional case 
    of system \eqref{1.1.main_system}, and present
    an appropriate numerical method to solve it. In their works the authors have used an abstract 
    scheme throughout Leray projection and classical techniques of optimal control theory. In 
    \cite{medjo2008optimal}, some theoretical and numerical aspects are presented for the optimal and 
    robust control of the Navier--Stokes equations. Additional information on optimal and robust control theory 
    for linear and nonlinear systems can be found in \cite{bewley1997optimal}, \cite{medjo2008optimal},
    \cite{green2012linear}, and references therein.
    
    In the context of hierarchical control (that is $\psi\equiv 0$ in \eqref{1.1.main_system}), some 
    recent works such as \cite{guillen2013approximate, araruna2014approximate, araruna2015stackelberg,
    limaco2009remarks} and \cite{hernandez2016robust} show a strategy with a leader and 
    follower controls for different equations. Some older results on a 
 	Stackelberg--Nash control strategy were proved by J. Diaz and J-L. Lions in 		
	\cite{diaz2004approximate} for a linear parabolic problem and by J. Limaco et al. 
  	for a linear parabolic problem with moving boundaries \cite{limaco2009remarks}. In both cases, the objective of the
  	leader control is an approximate controllability result. In the case of linear fluid
    models some approximate controllability of Stackelberg--Nash strategies  started with the result of F. Guill\'en--
    Gonz\'alez et. al for the Stokes system \cite{guillen2013approximate}, and were  extended by 
    F. Araruna et al. for linearized micropolar fluids \cite{araruna2014approximate}. In 
    \cite{araruna2014approximate} the main arguments are based on a  Fenchel-- 
    Rockafeller dual variational principle \cite{zalinescu2002convex}. For semilinear parabolic
    equations, a Stackelberg--Nash strategy with exact controllability for the leader control is proved in
     \cite{araruna2015stackelberg} using Carleman inequalities.
    
  	In our work, we follow the ideas introduced in \cite{hernandez2016robust} for the Navier--Stokes equations with
  	Dirichlet boundary conditions.  However the non linearity of \eqref{1.1.main_system} will allows only to obtain a
  	local null controllability result for the leader control.     
%
%

	Let us now introduce  the usual spaces in the context of incompressible fluids (\cite{Temam}): 
\begin{equation*}
    \begin{aligned}
    \begin{array}{lll}
    H:=&\{ u\in L^{2}(\Omega)^{N}: \nabla\cdot u=0, \,\,\text{in}\,\, \Omega,\,\,
    \, u\cdot n=0 \,\,\text{on}\,\, \partial\Omega \},\\
    V:=&\{ u\in H_{0}^{1}(\Omega)^{N}: \nabla\cdot u=0 \,\,\text{in}\,\, \Omega\}.
    \end{array}
    \end{aligned}
\end{equation*}
    
    Following the scheme for the robust control
    problem given in \cite{abergel1990some,bewley2000general}, the general space for the control
    functions and the disturbance $\psi$ in the right-hand side of
    \eqref{1.1.main_system} is $L^2(0,T;H)$.\\ 
    Now, we focus our attention on the control problem we are interested in.

\subsection{\normalsize The main problem}
    Given $h\in L^2(\omega\times (0,T))^N$ a (leader) control, we   consider the  \textit{secondary} cost functional
\begin{equation}\label{preli_functional_sec1}
    J_{r}(\psi,v;h):=\frac{\mu}{2}\iint\limits_{\mathcal{O}_d\times(0,T)}|y-y_{d}|^2         
    dxdt+\frac{1}{2}\Biggl(\ell^2\iint\limits_{\mathcal{O}\times(0,T)}\chi_{\mathcal{O}}|v|^2 dxdt
    -\gamma^2\iint\limits_{Q}|\psi|^2 dxdt \Biggr),
\end{equation}
    
   	where $\ell, \gamma,\mu>0$ are constants, $\mathcal{O}_d$ is an open subset of
    $\Omega$, which represents a observability domain,  and 
    $y_{d}\in L^2(0,T;L^2(\mathcal{O}_d)^N)$ is given. The constant $\mu $     
    arises from the physical parameters that govern the motion of fluids such     
    as  viscosity, characteristic length and characteristic velocity. The parameters    
    $\ell, \gamma$ are included to make the cost functional consistent and to 
    account the relative weight of each term. Note that the sign of the term associated to the
    disturbance is opposite to the sign used for the control, this is because we minimize with respect to the
    control $v$ meanwhile simultaneouly maximize with respect to the disturbance $\psi$. 
    From another perspective, the term $-\gamma^2\|\psi\|^2_{L^2(Q)^N}$ constrains the 
    magnitude of the disturbance function in the maximization with respect to $\psi$ and, the
    term associated to  	$\ell^2\|v\|^2_{L^2(Q)^N}$ constrains the magnitude of the control in the 
    minimization with respect to $v$.

    To explain the robust Stackelberg control problem, we will consider the following two 
    subproblems:
\begin{itemize}
    \item [i)] \textit{First problem.}  For every fixed leader control $h$, solve 
    the robust control problem for the nonlinear system   
    \eqref{1.1.main_system}, that is, find the best control $v$ in the presence of the 
    disturbance $\psi$ which maximally spoils the  follower control for the     
    Navier--Stokes system \eqref{1.1.main_system}. The robust control         
    problem to be solved is given in the following definition.

    \begin{Definition}\label{cp_def_nonlinearcase} Let 
    $h\in L^2(\omega\times(0,T))^N$ be fixed.
    The disturbance 
    $\overline{\psi}\in L^2(Q)^N$, the control  
    $\overline{v}\in L^2(Q)^N$, and the solution 
    $\overline{y}=y(h,\overline{v}(h),\overline{\psi}(h))$ of 
    \eqref{1.1.main_system} associated with 
    $(\overline{\psi}(h),\overline{v}(h))$
    are said to solve the robust control problem when a saddle point 
    $(\overline{\psi}(h),\overline{v}(h))$ of the cost functional defined in 
    \eqref{preli_functional_sec1} is reached, that is, if 
    \begin{equation}\label{preli_saddlecond_nonlinear_case}
    J_{r}(\psi,\overline{v}(h);h)\leq J_{r}(\overline{\psi}(h),\overline{v}(h);h)
    \leq  J_{r}(\overline{\psi}(h),v(h);h),\quad 
    \forall (\psi,v)\in L^2(Q)^{N\times N}.
    \end{equation}
    In this case, 
    $$J_{r}(\overline{\psi}(h),\overline{v}(h);h)=\max\limits_{\psi\in L^2(Q)^N}
    \min\limits_{v\in L^2(Q)^N} J_{r}(\psi,v;h)=\min\limits_{v\in L^2(Q)^N}\max    
    \limits_{\psi\in L^2(Q)^N} J_{r}(\psi,v;h).$$
    \end{Definition}

    \item [ii)] \textit{Second problem.}  Once the saddle point has been identified for each 
    leader control $h$, this is, once the existence of the  saddle point  $(\overline{\psi}(h),\overline{v}(h))$ for
    every leader control $h$ is guarantied, we deal with the problem of finding  the control $h$ of minimal norm
    satisfying  null controllability constraints. More precisely, we look for an optimal control 
    $\overline{h}$ such that 
    \begin{equation}\label{cp_second_subproblem}
    J(\overline{h})=\min\limits_{h}\frac{1}{2}\iint\limits_{\omega\times(0,T)}
    |h|^2 dxdt,\quad \mbox{subject to the restriction}\quad y(\cdot, T)=0\,\,\,\mbox{in}\,\,\Omega.
    \end{equation}
\end{itemize}

    Our main result on the robust hierarchic control is given in the following     
    theorem.
\begin{teo}\label{th_mainresult} Assume that 
    $\omega\cap \mathcal{O}_d\neq \emptyset$. Then, for every $T>0$ and 
    $\mathcal{O},\omega\subset \Omega$ open subsets such that $\mathcal{O}\cap \omega=\emptyset$, 
	there exist $\gamma_0,\ell_0,\delta$ and a positive function $\rho=\rho(t)$ blowing up $t=T$ such that 
	for any $\gamma\geq \gamma_0,\,\ell\geq\ell_0,\, y_0\in V$ and $y_d\in L^2(0,T;L^2(\mathcal{O}_d)^N)$
    satisfying
	\begin{equation}\label{hypothesis.weight.mainresult}
		\|y_0\|_{V}\leq \delta	\quad\mbox{and}\quad 
		\displaystyle\iint\limits_{\mathcal{O}_d\times(0,T)}\rho^2(t)|y_d|^2 dxdt<+\infty,
	\end{equation}	    
    we can find a leader control $h\in  L^2(0,T;L^2(\omega)^N)$
 	and an unique saddle point $(\overline{\psi},\overline{v})$ on $L^2(Q)^N\times         
    L^2(0,T;L^2(\mathcal{O})^N)$ and an associated solution $(y,p)$ to \eqref{1.1.main_system} verifying 
  	$y(\cdot,T)=0 \mbox{ in } \Omega.$
\end{teo}

    In order to prove Theorem \ref{th_mainresult}, we shall mainly consider two steps: a) the
    robust control results established in 
    \cite{bewley2000general} allow us to solve the mentioned-above \textit{first problem}. Here, as consequence of
    the nonlinearity given by the convection term, constrains either over small data or small time are 
    necessary in order to obtain the robust control; b) The hierarchical control (\textit{second problem}),
	where the main tools will be news Carleman estimates and fixed point arguments for solving the local null
	controllability associated to the leader control.    
      
    The rest of the paper is organized as follows. In Section 2, we present the general scheme of the robust 
    control problem for the system \eqref{1.1.main_system}. 
    In the first subsection we present the existence and characterization of the robust control for the linearized system
    (Stokes equation) and in the second subsection the same result for the nonlinear case. In section 3, we solve the
    robust Stackelberg strategy for the Stokes case. That is, we	
    prove the null controllability for the coupled Stokes system that arises as characterization of the 
    robust control problem.  In Section 4, we end the proof of Theorem \ref{th_mainresult} throughout an
    inverse function theorem of the Lyusternik's kind.

\section{\normalsize The robust control problem}
	As mentioned in the previous section, the main objetive in  robust control is to determine 
    the best control function
    $v\in L^2(\mathcal{O}\times(0,T))^N$ in the presence of the disturbance $\psi\in L^2(Q)^N$ which maximally 
    spoils the control. In this section we present some lemmas on the existence, uniqueness and 
    characterisation of a solution to the robust control problem established in Definition
    \ref{cp_def_nonlinearcase}.\\
    The proof of the existence of a solution $(\overline{\psi},\overline{v})$ to the robust control problem is
    based on the following result. The interested reader can see \cite{ekeland1999convex} for more details. 
\begin{lema}\label{2.existence.Ekeland_Temam}
    Let $\mathcal{J}$ be a functional defined on $X\times Y$, where $X$ and $Y$ are non--empty. closed, 
    unbounded convex sets. If $\mathcal{J}$ satisfies
    \begin{enumerate}
        \item [a)] $\forall \psi\in X,\, v\longmapsto \mathcal{J}(\psi,v)$ is convex lower semicontinuous.
        \item [b)] $\forall v\in Y,\, \psi\longmapsto \mathcal{J}(\psi,v)$ is concave upper semicontinuous.
        \item [c)] $\exists\psi_0\in X$\,\, such that $\lim\limits_{\|v\|_Y\to\infty}\mathcal{J}(\psi_0,v)
                    =+\infty$.
        \item [d)] $\exists v_0\in Y$\,\, such that $\lim\limits_{\|\psi\|_X\to\infty}\mathcal{J}(\psi,v_0)
                    =-\infty$.
    \end{enumerate}
   Then the functional $\mathcal{J}$ has a least one saddle point $(\overline{\psi},\overline{v})$ and
   $$\mathcal{J}(\overline{\psi},\overline{v})=\min\limits_{v\in Y}\sup\limits_{\psi\in X}\mathcal{J}(\psi,v)
   =\max\limits_{\psi\in X}\inf\limits_{v\in Y}\mathcal{J}(\psi,v).$$    
\end{lema}

\subsection{\normalsize Linear problem} 
	In this section we will treat the corresponding robust Stackelberg strategy for the linearized system. That is we
	will consider the Stokes system 
\begin{equation}\label{1.1.stokes_sys}
    \left\{
    \begin{aligned}
    \begin{array}{llll}
    y_{t}-\Delta y+\nabla p=h1_{\omega}+v\chi_{\mathcal{O}}+\psi  &
    \text{ in }& Q,\\    
    \nabla \cdot y=0&  \text{ in }& Q, \\
    y=0&\text{ on }&\Sigma, \\
    y(\cdot,0)=y_0(\cdot)& \text{ in }&\Omega.
    \end{array}
    \end{aligned}\right.
\end{equation}    
    
	We have the following result:
   	\begin{lema}\label{2.lemma.existence.saddlepoint}
    Let $h\in L^2(\omega\times(0,T))^N$ be fixed. There exists $\gamma_0>0$ such that for every $\gamma>\gamma_0$,
    there exists a saddle point 
    $(\overline{\psi},\overline{v})$ and the corresponding solution 
    $y(h,\overline{\psi},\overline{v})$ of \eqref{1.1.stokes_sys} such that
    $$J_{r}(\psi,\overline{v};h)\leq J_{r}(\overline{\psi},\overline{v};h)
    \leq J_{r}(\overline{\psi},v;h),\,\, \forall (\psi,v)\in L^2(Q)^N\times L^2(0,T;L^2(\mathcal{O})^N).$$
\end{lema}
	The proof of Lemma \ref{2.lemma.existence.saddlepoint} follows as in \cite{bewley2000general}  where the authors
	used Lemma \ref{2.existence.Ekeland_Temam} with $X=Y=L^2(Q)^N$  to prove the existence of a saddle point for 
	a slightly different cost functional $\mathcal{J}$. 
    As consequence of this result, the existence of a solution $(\overline{\psi},
     \overline{v})$ to our robust control problem is guaranteed.  
\begin{Obs}
 	In lemma \ref{2.lemma.existence.saddlepoint}, if the condition on $\gamma$ is 
    not met, we cannot prove the existence of the saddle point.
    On the other hand, it is known that the existence of a saddle point for the functional $J_{r}$ implies
    that for any $\psi\in L^2(Q)^N,\, v\in L^2(0,T;L^2(\mathcal{O})^N)$  	
	$$\frac{\partial J_{r}}{\partial \psi}(\overline{\psi},\overline{v})\cdot\psi=0,\quad 
	\frac{\partial J_{r}}{\partial v}(\overline{\psi},\overline{v})\cdot v=0,$$
	where
	\begin{equation}
		\frac{\partial J_{r}}{\partial \psi}(\overline{\psi},\overline{v})\cdot {\psi}
		=\iint\limits_{\mathcal{O}_d\times(0,T)}(y-y_d){w}_\psi dxdt
		-\mathcal{\gamma}^{2}\iint\limits_{\mathcal{O}\times(0,T)}\psi\overline{\psi}dxdt
	\end{equation}
	and 
	\begin{equation}
		\frac{\partial J_{r}}{\partial v}(\overline{\psi},\overline{v})\cdot {v}
		=\iint\limits_{\mathcal{O}_d\times(0,T)}(y-y_d){w}_v dxdt
		+\mathcal{\ell}^{2}\iint\limits_{\mathcal{O}\times(0,T)}\chi_{\mathcal{O}}v\overline{v}dxdt,
	\end{equation}
	and ${w}_\psi,\, {w}_v$ are the G\^ateaux derivatives of $y$ solution to \eqref{1.1.stokes_sys} in 
	the directions $\psi$ and $v$ respectively.
\end{Obs}
	Finally, in order to characterize the robust control problem, we introduce the linear
	 adjoint system \textcolor{red}{to} (\ref{1.1.stokes_sys}) with right--hand side related with $J_r$, that is, 
	 we consider   
\begin{equation}
	\left\{
	\begin{aligned}
	\begin{array}{llll}
		-z_{t}-\Delta z+\nabla \pi_z
		        =\mu(y-y_d)\chi_{\mathcal{O}_d} &\text{ in }& Q,\\
		        \nabla \cdot z=0  &\text{ in }& Q, \\
		        z=0&\text{ on }&\Sigma, \\
		       z(\cdot,T)=0& \text{ in }&\Omega.
	\end{array}
	\end{aligned}\right.
\end{equation}	
	In the following result we characterize the saddle point $(\overline{v},\overline{\psi})$ in terms of $z$.
	The interested reader can consult \cite{bewley2000existence} for more details. 
	
\begin{lema}\label{2.teo.robust_charac.saddle_point1}
    Let $h\in L^2(\omega\times(0,T))^N$ and $y_0\in V$ be given. Suppose that 
    $(\overline{\psi},\overline{v})$ is the solution to the robust control 
    problem stated in Definition \ref{cp_def_nonlinearcase}. Then  
    \begin{equation*}
        \overline{\psi}=\frac{1}{\gamma^2}z\quad\mbox{and}\quad 
        \overline{v}=-\frac{1}{\ell^2}z\chi_{\mathcal{O}},
    \end{equation*}
    where  $\gamma$ is sufficiently large and  the pair $(y,z)$ solves the following coupled system:
    \begin{equation}\label{robust_coupled_system_linear}
        \left\{
        \begin{aligned}
        \begin{array}{llll}
        y_{t}-\Delta y+\nabla \pi_y
        =h1_{\omega}+(-\ell^{-2}\chi_{\mathcal{O}}            
        +\gamma^{-2})z &\text{ in }& Q,\\
        -z_{t}-\Delta z+\nabla \pi_z
        =\mu(y-y_d)\chi_{\mathcal{O}_d} &\text{ in }& Q,\\
        \nabla \cdot y=0, \nabla \cdot z=0  &\text{ in }& Q, \\
        y=z=0&\text{ on }&\Sigma, \\
        y(\cdot,0)=y_0(\cdot),\quad z(\cdot,T)=0& \text{ in }&\Omega.
        \end{array}
        \end{aligned}\right.
     \end{equation}
\end{lema}

\subsection{\normalsize Nonlinear problem}

    The analysis is similar to the previous one for the linear case. However, 
    it is well known that the theory of the 
    Navier--Stokes equations is complete in two--dimensional spaces, which do not occur in three--dimensional
    spaces. Roughly speaking, in three dimensions, the existence of a robust control is restricted to cases of 
    either small data or small $T$. Additionally, the nonlinearity will require new assumptions on the parameter $\ell$. 
    Under the constraint of small data, we need to impose
    the following condition: there exists $\delta>0$ such that,
    for every $(v\chi_{\mathcal{O}},\psi)\in L^2(Q)^{N\times N}$ and $y_0\in V$ 
\begin{equation}\label{robuscontrol.smalldata}
    \|v\chi_{\mathcal{O}}\|_{L^2(Q)^N}+\|\psi\|_{L^2(Q)^N}\leq \delta \quad \mbox{and}\quad \|y_0\|_V\leq \delta
\end{equation}     
    holds.
    \newpage 
\begin{lema}\label{lemma.nonlinearcase_RC}
	Let $h\in L^2(\omega\times(0,T))^N$ be fixed.
    \begin{enumerate} 
        \item [i)] Case $N=2$. There exist constants $\gamma_0>0$ and $\ell_0>0$ 
        such that for every 
        $\gamma>\gamma_0$ and $\ell>\ell_0$, there exists 
        $(\overline{\psi},\overline{v})$ on $L^2(Q)^N\times L^2(0,T;L^2(\mathcal{O})^N)$ and 
     	the  associated solution to \eqref{1.1.main_system}
        $y=y(h,\overline{v},\overline{\psi})$ such that 
        $$J_{r}(\psi,\overline{v};h)\leq J_{r}(\overline{\psi},\overline{v};h)
    	\leq J_{r}(\overline{\psi},v;h),\,\, \forall\, (\psi,v)\in 
    	L^2(Q)^N\times L^2(0,T;L^2(\mathcal{O})^N).$$ That is,  $(\overline{\psi},\overline{v})$ is a saddle point 
    	of $J_r$.
    	
        \item [ii)] Case $N=3$. Under the hypothesis of the case $N=2$,  
        and that either $y_0\in V$ and  $(v\chi_{\mathcal{O}},\psi)\in L^2(Q)^{N\times N}$
        satisfies \eqref{robuscontrol.smalldata}, or that $t=T$ is small, then there exists  
		$(\overline{\psi},\overline{v})\in L^2(Q)^N\times L^2(0,T;L^2(\mathcal{O})^N)$ a saddle point of $J_{r}$.
    \end{enumerate}
\end{lema}

 	Analogously to the linear case,  we give the characterization of the robust control problem in the following
 	result.
\begin{lema}\label{2.teo.robust_charac.saddle_point2}
    Let $h\in L^2(\omega\times(0,T))$ and $y_0\in V$ be given. Then, there exist positive constants
    $\gamma_0,\, \ell_0$
    such that if $\gamma>\gamma_0,\, \ell>\ell_0$, the solution $(\overline{v},\overline{\psi})$
    to the robust control problem stated in Definition \ref{cp_def_nonlinearcase} exists and is unique. Furthermore,  
    $(\overline{v},\overline{\psi})$ is 
    characterized by
    \begin{equation}\label{oooooo}
        \overline{\psi}=\frac{1}{\gamma^2}z\quad\mbox{and}\quad 
        \overline{v}=-\frac{1}{\ell^2}z\chi_{\mathcal{O}},
    \end{equation}
    where  $z$ is the second component of $(y,z)$ solution to the
    following coupled system:
    \begin{equation}\label{robust_coupled_system_nonlinear}
        \left\{
        \begin{aligned}
        \begin{array}{llll}
        y_{t}-\Delta y+(y\cdot\nabla)y+\nabla \pi_y
        =h1_{\omega}+(-\ell^{-2}\chi_{\mathcal{O}}            
        +\gamma^{-2})z &\text{ in }& Q,\\
        -z_{t}-\Delta z+(z\cdot\nabla^t)y-(y\cdot\nabla)z+\nabla \pi_z
        =\mu(y-y_d)\chi_{\mathcal{O}_d} &\text{ in }& Q,\\
        \nabla \cdot y=0, \nabla \cdot z=0  &\text{ in }& Q, \\
        y=z=0&\text{ on }&\Sigma, \\
        y(\cdot,0)=y_0(\cdot),\quad z(\cdot,T)=0& \text{ in }&\Omega.
        \end{array}
        \end{aligned}\right.
     \end{equation}
\end{lema}
	The proof of Lemma \ref{lemma.nonlinearcase_RC} and Lemma  \ref{2.teo.robust_charac.saddle_point2} 
	can be found in \cite{bewley2000general}.
\section{\normalsize Controllability}
	In the previous sections we saw that the robust control is characterized in such a way that a coupled system needs 
	to be solved. In order to establish  a Stackelberg strategy requiring the leader control to drive the equation to
	zero we need to find $h\in L^2(\omega\times (0,T))^N$ such that the corresponding $y$ solution to
	\eqref{robust_coupled_system_linear} (in the linear case) or to \eqref{robust_coupled_system_nonlinear} 
	(in the nonlinear case), satisfies $y(T)=0$. To achieve this objectives, we will obtain first the result in the
	linear case. To this aim we will prove an observability inequality for the adjoint system to
	\eqref{robust_coupled_system_linear} by means of Carleman estimates. The nonlinear case will be obtained by a fixed
	point argument. The next subsection will be devoted to the obtention of the Carleman inequalities. 
\subsection{\normalsize Carleman inequalities}

     We first define several weight functions which will be useful in the sequel.  
     Let  $\omega_0$ be a nonempty open subset of $\mathbb{R}^N$ such that
     $\omega_0\Subset{\omega}\cap \mathcal{O}_d$ 
     and  $\eta\in C^2(\overline{\Omega})$ such that
     $$|\nabla \eta|>0 \mbox{ in }\overline{\Omega}\setminus\omega_0,\,\,\,\, \eta>0
     \mbox{ in }\Omega\,\,\, \mbox{ and }\,\, \eta \equiv 0 \mbox{ on }\partial
     \Omega.$$ 
     The existence of such a function $\eta$ is proved in \cite{Fursikov}. Then, for some positive 
     real number $\lambda$, we consider the following weight functions:
\begin{equation}\label{carleman_weights}
    \begin{aligned}
    &\alpha(x,t) = \dfrac{e^{12\lambda\|\eta\|_{\infty}}-e^{\lambda(10\|\eta\|_{\infty}+\eta(x))}} 
    {(t(T-t))^{5}},\quad \xi(x,t)=\dfrac{e^{\lambda(10\|\eta\|_{\infty}+\eta(x))}}{(t(T-t))^{5}},\\
    &\alpha^*(t) = \max_{x\in\overline{\Omega}} \alpha(x,t),\quad \quad
     \xi^*(t) = \min_{x\in\overline{\Omega}} \xi(x,t),\\
    &\widehat\alpha(t) = \min_{x\in\overline{\Omega}} \alpha(x,t),\quad 
    \quad \,\,\,\,\widehat\xi(t) = \max_{x\in\overline{\Omega}} \xi(x,t).
    \end{aligned}
\end{equation}
	These weight functions have been used by M. Gueye in \cite{gueye2013insensitizing} and S. Guerrero in 
    \cite{guerrero2007controllability} to obtain
    Carleman estimates for a Stokes coupled system similar to the presented in our work.\\ 
       
   We consider now the non homogeneous adjoint system to \eqref{robust_coupled_system_linear}:
\begin{equation}\label{carleman_adjoint_system}
    \left\{
    \begin{aligned}
    \begin{array}{llll}
        -\varphi_{t}-\Delta \varphi+\nabla \pi_{\varphi}
        =g_1+\mu\theta\chi_{\mathcal{O}_d} &\text{ in }& Q,
        \\
        \theta_{t}-\Delta \theta+\nabla \pi_{\theta}
        =g_2 -\ell^{-2}\varphi \chi_{\mathcal{O}}+\gamma^{-2}\varphi &\text{ in }& Q,
        \\
        \nabla \cdot  \varphi=0,\, \nabla \cdot \theta=0  &\text{ in }& Q, \\
        \varphi=\theta=0&\text{ on }&\Sigma, \\
        \varphi(\cdot,T)= \varphi_T(\cdot),\, \theta(\cdot,0)=0& \text{ in }&\Omega,
    \end{array}
    \end{aligned}
    \right.
\end{equation}
    where $g_1,g_2\in L^2(Q)^N$ and $\varphi_T\in H$.\\
     
    Our Carleman estimate is given in the following proposition. In what follows, the constants $a_0$ and
    $m_0$ are fixed, and satisfy  
\begin{equation}\label{constantes}
		\frac{5}{4}\leq a_0<a_0+1<m_0<2a_0,\quad m_0<2+a_0. 
\end{equation}

\begin{proposition}\label{carleman_carleman1}
    Assume that $\omega\cap \mathcal{O}_d\neq \emptyset$ and that $\ell$ and $\gamma$ are large enough.
    Then, there exist a constant  $\overline{\lambda}$ such that for any $\lambda\geq\overline{\lambda}$
    exist two constants $\overline{s}(\lambda)>0$ and $C=C(\lambda)>0$ depending only on 
    $\Omega$ and $\omega$  such that  for any $g_1,g_2\in L^2(Q)^N$ and any $\varphi_T\in H$,  the solution of
    \eqref{carleman_adjoint_system} satisfies
    \begin{equation}\label{carleman_obs_ineq1}
    \begin{array}{lll}
        &\displaystyle\iint\limits_{Q}e^{-2s\alpha-2a_0s\alpha^*}(s\lambda^2\xi|\nabla(\nabla\times\theta)|^2
        +s^3\lambda^4\xi^3|\nabla\times\theta|^2)dxdt\\
        &\hspace{1cm}
        +\displaystyle\iint\limits_{Q}e^{-2sm_0\alpha}(s\lambda^2\xi|\nabla\varphi|^2+s^3\lambda^4\xi^3|
        \varphi|^2+(s\xi)^{-1}|\Delta\varphi|^2)dxdt\\
        &\leq C\Biggl(s^{15}\lambda^{24}\displaystyle\iint\limits_{\omega\times(0,T)}
        e^{-4a_0s\alpha^*+2(m_0-2)s\alpha^*}(\hat\xi)^{15}|\varphi|^2 dxdt\\
        &\hspace{2cm}+s^{5}\lambda^{6}\displaystyle\iint\limits_{Q}
        e^{-2s\hat\alpha-2a_0s\alpha^*}(\hat\xi)^{5}|g_1|^2dxdt
        +\displaystyle\iint\limits_{Q}e^{-2a_0s\alpha^*}|g_2|^2dxdt\Biggr),
    \end{array}
    \end{equation}
    for any $s\geq \overline{s}$.
\end{proposition}

    Before giving the proof of Proposition \ref{carleman_carleman1}, we recall some technical
    results. We first present a Carleman inequality proved in \cite{FC-B-G-P} for a general  heat equation
    with Fourier boundary conditions. Let us introduce the system
\begin{equation}\label{carleman_FCBGP06}
    \left\{
    \begin{aligned}
    \begin{array}{llll}
    -u_{t}-\Delta u=f_1+\nabla\cdot f_2  &\text{ in }& Q,\\
    (\nabla u+f_2)\cdot n=f_3&\text{ on }&\Sigma, \\
    u(\cdot,T)=u_{T}(\cdot)& \text{ in }&\Omega,
    \end{array}
    \end{aligned}\right.
\end{equation}
    where $f_1\in L^2(Q), f_2\in L^2(Q)^N$ and $f_3\in L^2(\Sigma)$. We have:
\begin{lema}\label{carleman_lemaFCBGP06}
    Under the previous assumptions on $f_1,f_2$ and $f_3$, there exist positive constants
    $\overline{\lambda}, \sigma_1,\sigma_2$ and $C$, only depending on $\Omega$ 
    and $\omega$, such that, for any $\lambda\geq\overline{\lambda}$, any 
    $s\geq \overline{s}=\sigma_1(e^{\sigma_2\lambda}T+T^2)$ and any 
    $u_T\in L^2(\Omega)$, the weak solution to \eqref{carleman_FCBGP06} satisfies
    \begin{equation}\label{carlemanFCBGP06}
        \begin{aligned}
        &\iint\limits_{Q}e^{-2s\alpha}[s^3\lambda^4\xi^3|u|^2 +s\lambda^2\xi|\nabla u|^2]dxdt
        \leq C\Bigl(\iint\limits_{Q}e^{-2s\alpha}(|f_1|^2+s^2\lambda^2\xi^2|f_2|
        ^2)dxdt\\
        &\hspace{3cm}+s\lambda\iint\limits_{\Sigma}e^{-2s\alpha}\xi|f_3|^2d\sigma dt 
        +s^3\lambda^4\iint\limits_{\omega_0\times (0,T)}e^{-2s\alpha}\xi^3|u|^2dxdt\Bigr).
        \end{aligned}
    \end{equation}
\end{lema}
    The second result holds for the solutions of a Stokes system with Dirichlet boundary 
    conditions. The interested reader can see \cite{FCGIP04} for more details.
\begin{lema}\label{carleman_lemaFCGOIP04}
    Let $u_0\in V$ and  $f_4\in L^2(Q)^N$. Then, there exists a constant $C(\Omega,\omega,T)>0$ 
    such that the solution $(u,p)\in L^2(0,T;H^2(\Omega)^N\cap V)\cap L^{\infty}(0,T;V)\times 
    L^2(0,T;H^1(\Omega))$, with $\displaystyle\int\limits_{\omega_0}p(x,t)dx=0$, of
    \begin{equation*}
        \left\{
        \begin{aligned}
        \begin{array}{llll}
        u_{t}-\Delta u+\nabla p=f_4  &\text{ in }& Q,\\
        \nabla\cdot u=0 &\text{ in }& Q,\\
         u=0&\text{ on }&\Sigma, \\
        u(\cdot,0)=u_{0}(\cdot)& \text{ in }&\Omega,
        \end{array}
        \end{aligned}\right.
    \end{equation*} 
    satisfies 
    \begin{equation}\label{carleman_lemaFCGIOP04_inequality}
        \begin{array}{lll}
        &\displaystyle\iint\limits_{Q}e^{-2s\alpha}(s\lambda^2\xi|\nabla u|^2+s^3\lambda^4\xi^3|u|^2)dxdt\\
        &\leq C\Bigl( s^{16}\lambda^{40}\displaystyle\iint\limits_{\omega\times(0,T)}e^{-8s\hat{\alpha}
        +6s\alpha^*}(\hat{\xi})^{16}|u|^2 dxdt
        +s^{15/2}\lambda^{20}\displaystyle\iint\limits_{Q}e^{-4s\hat{\alpha}+2s\alpha^*}(\hat{\xi})^{15/2}|
        f_4|^2 dxdt\Bigr),
        \end{array}
    \end{equation}
    for any $\lambda\geq C$ and $s\geq C(T^5+T^{10})$.
\end{lema}
   
\begin{Obs}
    In \cite{FC-B-G-P} and \cite{FCGIP04} slightly different weight functions are used to prove the above     
    results. However, the inequality remains valid since the key point of the proof is that $\alpha$ 
    goes to $0$ when $t$ tends to $0$ and $T$.
\end{Obs}
    The next result concerns the regularity of the solutions to the Stokes system, see 
    \cite{Temam} and \cite{guerrero2007controllability} for more details.

\begin{lema}\label{carleman_regularityH3}
    Let $a\in\mathbb{R}$ and $B\in\mathbb{R}^N$ be constant and let $f_5\in L^2(0,T;V)$. Then, there exists a 
    unique solution 
    $$u\in L^2(0,T;H^3(\Omega)^N\cap V)\cap H^1(0,T;V)$$ 
    for the Stokes system 
    \begin{equation}\label{system.regularityH3}
        \left\{
        \begin{aligned}
        \begin{array}{llll}
        u_{t}-\Delta u+au+B\cdot\nabla u+\nabla p=f_5  &\text{ in }& Q,\\
        \nabla\cdot u=0 &\text{ in }& Q,\\
         u=0&\text{ on }&\Sigma, \\
        u(\cdot,0)=0& \text{ in }&\Omega,
        \end{array}
        \end{aligned}\right.
    \end{equation} 
    for some $p\in L^2(0,T;H^2(\Omega))$, and there exists a constant $C>0$ such that
    \begin{equation}\label{carleman_regularityH3_ineq}
        \|u\|_{L^2(0,T;H^3(\Omega)^N)}+\|u\|_{H^1(0,T;L^2(\Omega)^N)}\leq C\|f_5\|_{L^2(0,T;H^1(\Omega)^N)}.
    \end{equation}
    Moreover, if we assume that $a\equiv B\equiv 0$ and $f_5\in L^2(Q)^N$, $u$ is actually, together a 
    pressure $p$, the strong solution of \eqref{system.regularityH3}, i.e., 
    $(u,p)\in L^2(0,T;H^2(\Omega)^N)\cap L^{\infty}(0,T;V)\cap H^1(0,T;H)\times L^2(0,T;H^1(\Omega))$. 	
    Furthermore, there exists a constant $C>0$ such that 
    \begin{equation}\label{carleman_regularityH2_ineq}
        \|u\|_{L^2(0,T;H^2(\Omega)^N)}+\|u\|_{L^{\infty}(0,T;V)}+\|u\|_{H^1(0,T;L^2(\Omega)^N)}
        \leq C\|f_5\|_{L^2(Q)^N}.
    \end{equation}    
\end{lema}
     Now, we give the proof of Proposition \ref{carleman_carleman1}.
\subsection{\normalsize Proof of Proposition \ref{carleman_carleman1}}     
{\bf{Carleman estimate for $\theta$}}\\

    Let define $\theta^*:=\rho^*\theta,\quad \pi^*:=\rho^*\pi$, where $\rho^*=\rho^*(t)=e^{-a_{0}s\alpha^*}$ 
    and  $a_0$ fixed satisfying \eqref{constantes}.
    From \eqref{carleman_adjoint_system}, $(\theta^*,\pi^*)$ is the solution of the following system
\begin{equation*}
    \left\{
        \begin{aligned}
        \begin{array}{llll}
        \theta^*_{t}-\Delta \theta^*+\nabla \pi^*
        =\rho^*g_2+
        \rho^*(-\ell^{-2}\varphi\chi_{\mathcal{O}}+\gamma^{-2}\varphi)+\rho^*_{t}\theta &\text{ in }& Q,\\
        \nabla\cdot \theta^*=0 &\text{ in }& Q,\\
         \theta^*=0&\text{ on }&\Sigma, \\
        \theta^*(\cdot,0)=0& \text{ in }&\Omega.
        \end{array}
    \end{aligned}\right.
\end{equation*} 
    Now, we decompose $(\theta^*,\pi^*)$ as follows: 
    \begin{equation}\label{carleman_decompose}
        (\theta^*,\pi^*)=(\hat{\theta},\hat{\pi})+(\tilde{\theta},\tilde{\pi}),
    \end{equation}    
    where $(\hat{\theta},\hat{\pi})$ and $(\tilde{\theta},\tilde{\pi})$ solve respectively
\begin{equation}\label{carleman_sys_thetatilde}
    \left\{
        \begin{aligned}
        \begin{array}{llll}
        \tilde{\theta}_{t}-\Delta\tilde{\theta}+\nabla \tilde{\pi}
        =\rho^*g_2+\rho^*(-\ell^{-2}\varphi\chi_{\mathcal{O}}+\gamma^{-2}\varphi)&\text{ in }& Q,\\
        \nabla\cdot \tilde{\theta}=0 &\text{ in }& Q,\\
        \tilde{\theta}=0&\text{ on }&\Sigma, \\
        \tilde{\theta}(\cdot,0)=0& \text{ in }&\Omega,
        \end{array}
        \end{aligned}
    \right.
\end{equation} 
    and
\begin{equation}\label{carleman_sys_thetahat}
    \left\{
        \begin{aligned}
        \begin{array}{llll}
        \hat{\theta}_{t}-\Delta\hat{\theta}+\nabla \hat{\pi}
        =\rho^*_{t}\theta&\text{ in }& Q,\\
        \nabla\cdot \hat{\theta}=0 &\text{ in }& Q,\\
        \hat{\theta}=0&\text{ on }&\Sigma, \\
        \hat{\theta}(\cdot,0)=0& \text{ in }&\Omega.
        \end{array}
        \end{aligned}
    \right.
\end{equation}    
    For system \eqref{carleman_sys_thetatilde} we will use Lemma \ref{carleman_regularityH3} and the 
    regularity result estimate \eqref{carleman_regularityH2_ineq}, meanwhile for the system 
    \eqref{carleman_sys_thetahat} we will use the ideas of both works \cite{guerrero2007controllability} and 
    \cite{gueye2013insensitizing}.\\
    We apply the operator $\nabla\times\cdot$ to the Stokes system satisfied by $\hat{\theta}$. Then, we have
    $$(\nabla\times\hat\theta)_{t}-\Delta(\nabla\times\hat\theta)
    =\rho^*_{t}(\nabla\times\theta)\quad\mbox{in}\,\,Q.$$
    Using Lemma \cite{FC-B-G-P} with $f_1=\rho^*_{t}(\nabla\times\theta)$, there exists a constant 
    $C=C(\Omega,\omega_0)>0$ such that 
\begin{equation}\label{carleman_estimateproof1}
    \begin{aligned}
         &\displaystyle\iint\limits_{Q}e^{-2s\alpha}(s\lambda^2\xi|\nabla(\nabla\times\hat\theta)|^2
        +s^3\lambda^4\xi^3|\nabla\times\hat\theta|^2)dxdt\\
         &\leq C\Bigl(\iint\limits_{Q}e^{-2s\alpha}|\rho^*_t|^2|\nabla\times\theta|^2dxdt\\
         &\hspace{2cm}+s\lambda\iint\limits_{\Sigma}e^{-2s\alpha}\xi\Bigl|
        \frac{\partial(\nabla\times\hat\theta)}{\partial n}\Bigr|^2d\sigma dt 
        +s^3\lambda^4\iint\limits_{\omega_0\times (0,T)}e^{-2s\alpha}\xi^3|\nabla\times\hat\theta|^2dxdt
        \Bigr),
    \end{aligned}
\end{equation}
    for any $\lambda\geq C$ and $s\geq C(T^{10}+T^9)$.\\
    Now, using the inequality $(a-b)^2\geq \frac{a^2}{2}-b^2$, for every $a,b\in \mathbb{R}$ with $a=\theta^*$
    and $b=\tilde{\theta}$, we get (recall that $\hat{\theta}=\theta^*-\tilde{\theta}$):
\begin{equation}\label{carleman_ine_aux1}
    \begin{array}{lll}
    &\displaystyle\frac{1}{2}\iint\limits_{Q}e^{-2s\alpha-2a_0s\alpha^*}
    (s\lambda^2\xi|\nabla(\nabla\times\theta)|^2+s^3\lambda^4\xi^3|\nabla\times\theta|^2)dxdt\\
    &\hspace{1cm}-\displaystyle\iint\limits_{Q}e^{-2s\alpha}
    (s\lambda^2\xi|\nabla(\nabla\times\tilde\theta)|^2+s^3\lambda^4\xi^3|\nabla\times\tilde\theta|^2)dxdt\\
    &\leq\displaystyle\iint\limits_{Q}e^{-2s\alpha}(s\lambda^2\xi|\nabla(\nabla\times\hat\theta)|^2
    +s^3\lambda^4\xi^3|\nabla\times\hat\theta|^2)dxdt.   
    \end{array}
\end{equation}
    The fact that $s^3\lambda^4e^{-2s\alpha}\xi^3$ and $s\lambda^2e^{-2s\alpha}\xi$ are upper bounded allow us 
    to estimate the terms associated to $|\nabla(\nabla\times\tilde\theta)|^2$ and 
    $|\nabla\times\tilde\theta|^2$ through 
    \eqref{carleman_regularityH2_ineq}. More precisely, we have:
\begin{equation}\label{carleman_aux2}
    \begin{array}{l}
   \displaystyle s^3\lambda^4\displaystyle\iint\limits_{Q}e^{-2s\alpha}\xi^3|\nabla\times\tilde\theta|^2dxdt
    + s\lambda^2\displaystyle\iint\limits_{Q}e^{-2s\alpha}\xi|\nabla(\nabla\times\tilde\theta)|^2dxdt\\
        \noalign{\smallskip}
  \leq \displaystyle C_{s,\lambda}\|\tilde\theta\|^2_{L^2(0,T;H^1(\Omega)^N)\cap L^2(0,T;H^2(\Omega)^N)}\\
      \noalign{\smallskip}
 \displaystyle
    \leq C_{s,\lambda}\|\rho^*g_2\|^2_{L^2(Q)^N}
    +C_{s,\lambda}\|\rho^*(-\ell^{-2}\varphi\chi_{\mathcal{O}}+\gamma^{-2}\varphi)\|^2_{L^2(Q)^N}, 
    \end{array}
\end{equation}    
    where $C_{s,\lambda}$ is a positive constant depending on $s$ and $\lambda$.\\
    
    On the other hand, taking into account that $|\rho^*_{t}|\leq CsT\rho^*(\xi^*)^{6/5}$ for every $s\geq C$,
    it follows that
    $$\iint\limits_{Q}e^{-2s\alpha}|\rho^*_t|^2|\nabla\times\theta|^2dxdt
    \leq Cs^2T^2\iint\limits_{Q}e^{-2s\alpha-2a_0s\alpha^*}(\xi^*)^{12/5}|\nabla\times\theta|^2dxdt,$$
    which can be absorbed by the first term in the right--hand side of \eqref{carleman_ine_aux1}, for every 
    $\lambda\geq 1,\, s\geq C$.\\
    From the identity $\theta^*=\hat{\theta}+\tilde{\theta}$ (recall \eqref{carleman_decompose}) and 
    \eqref{carleman_aux2}, it is easy to estimate the local term that appear in the right--hand side 
    of \eqref{carleman_estimateproof1} by:
\begin{equation}\label{carleman_aux4}
    \begin{aligned}
        &s^3\lambda^4\iint\limits_{\omega_0\times(0,T)}e^{-2s\alpha}\xi^3|\nabla\times\hat\theta|^2dxdt\\
        &\leq C\Bigl(s^3\lambda^4\iint\limits_{\omega_0\times(0,T)}e^{-2s\alpha}\xi^3
        (|\nabla\times\theta^*|^2+|\nabla\times\tilde\theta|^2)dxdt)\\
        &\leq         Cs^3\lambda^4\!\!\!\!
		\iint\limits_{\omega_0\times(0,T)}e^{-2s\alpha}\xi^3|\nabla\times\theta^*|^2
        dxdt+C_{s,\lambda}\|\rho^*g_2\|^2_{L^2(Q)^N}
        +C_{s,\lambda}\|\rho^*(-\ell^{-2}\varphi\chi_{\mathcal{O}}+\gamma^{-2}\varphi)\|^2_{L^2(Q)^N}.    
    \end{aligned}
\end{equation}
    Putting together \eqref{carleman_estimateproof1}, \eqref{carleman_ine_aux1} and \eqref{carleman_aux4}, 
    we have for the moment
\begin{equation}\label{carleman_end_step1}
    \begin{aligned}
        & \displaystyle\iint\limits_{Q}e^{-2s\alpha-2a_0s\alpha^*}
        (s\lambda^2\xi|\nabla(\nabla\times\theta)|^2+s^3\lambda^4\xi^3|\nabla\times\theta|^2)dxdt\\
        &\leq  
        C\Biggl(s^3\lambda^4\iint\limits_{\omega_0\times (0,T)}e^{-2s\alpha-2a_0s\alpha^*}\xi^3|
        \nabla\times\theta|^2dxdt+s\lambda\iint\limits_{\Sigma}e^{-2s\alpha}\xi\Bigl|
        \frac{\partial(\nabla\times\hat\theta)}{\partial n}\Bigr|^2d\sigma dt\Biggr)\\ 
        &\hspace{1cm}+C_{s,\lambda}\|\rho^*g_2\|^2_{L^2(Q)^N}
        +C_{s,\lambda}\|\rho^*(-\ell^{-2}\varphi\chi_{\mathcal{O}}+\gamma^{-2}\varphi)\|^2_{L^2(Q)^N},         
    \end{aligned}
\end{equation}
    for every $s\geq C$ and $\lambda\geq C$.\\
    
    The last step will be to estimate the boundary term 
    $$s\lambda\iint\limits_{\Sigma}e^{-2s\alpha}\xi\Bigl|
    \frac{\partial(\nabla\times\hat\theta)}{\partial n}\Bigr|^2d\sigma dt.$$
    To this end we follow the arguments of \cite{gueye2013insensitizing}. For brevity we omitted the calculus but 
    refer to \cite{gueye2013insensitizing}.\\
    Therefore, there exist $C_{s,\lambda}>0$ and $C>0$ such that
\begin{equation*}
    \begin{array}{lll}
    s\lambda\displaystyle\iint\limits_{\Sigma}e^{-2s\alpha}\xi\Bigl|\frac{\partial(\nabla\times\hat\theta)}
    {\partial n}\Bigr|^2d\sigma dt
    &\leq C\|\rho^*g_2\|^2_{L^2(Q)^N}
    +C_{s,\lambda}\|\rho^*(-\ell^{-2}\varphi\chi_{\mathcal{O}}+\gamma^{-2}\varphi)\|^2_{L^2(Q)^N}   
    \\
    &+\varepsilon \Biggl(\displaystyle\iint\limits_{Q}e^{-2s\alpha-2a_0s\alpha^*}
    (s\lambda^2\xi|\nabla(\nabla\times\theta)|^2+s^3\lambda^4\xi^3|\nabla\times\theta|^2)dxdt \Biggr),
    \end{array}
\end{equation*}
    for every $\varepsilon >0$.\\ 
    From the previous inequality and \eqref{carleman_end_step1} we conclude the following Carleman
    estimate for $\theta$:
\begin{equation}\label{c.estimatetheta.end}
    \begin{array}{lll}
    &\displaystyle\iint\limits_{Q}e^{-2s\alpha-2a_0s\alpha^*}
    (s\lambda^2\xi|\nabla(\nabla\times\theta)|^2+s^3\lambda^4\xi^3|\nabla\times\theta|^2)dxdt
    \leq  C\|\rho^*g_2\|^2_{L^2(Q)^N}\\
    &+C_{s,\lambda}\|\rho^*(-\ell^{-2}\varphi\chi_{\mathcal{O}}
    +\gamma^{-2}\varphi)\|^2_{L^2(Q)^N}
    +C s^3\lambda^4\displaystyle\iint\limits_{\omega_0\times (0,T)}e^{-2s\alpha-2a_0s\alpha^*}\xi^3|\nabla\times
    \theta|^2dxdt,
    \end{array}
\end{equation}
    for every  $s\geq C(T^5+T^{10})$ and $\lambda\geq C$.\\
     
\noindent{\bf{Carleman estimate for $\varphi$}}\\
    
    First, assuming that $\theta$ is given, we look at $\varphi$ as the solution of
\begin{equation*}
    \left\{
    \begin{aligned}
    \begin{array}{llll}
        -\varphi_{t}-\Delta \varphi+\nabla \pi_{\varphi}
        =g_1+\mu\theta\chi_{\mathcal{O}_d} &\text{ in }& Q,\\
        \nabla \cdot  \varphi=0 &\text{ in }& Q, \\
        \varphi=0&\text{ on }&\Sigma, \\
        \varphi(\cdot,T)=\varphi_T(\cdot)& \text{ in }&\Omega.
    \end{array}
    \end{aligned}
    \right.
\end{equation*}
    Now, we choose $\pi_{\varphi}$ such that $\displaystyle\int\limits_{\omega_0}\pi_{\varphi}dx=0$ 
    and we apply Lemma \ref{carleman_lemaFCGOIP04} with $f_4=g_1+\mu\theta\chi_{\mathcal{O}_d}$ and use the 
    weight function $m_0\alpha$ (instead of $\alpha$), where $a_0+1<m_0\leq 2a_0$ and $m_0\leq 2+a_0$.\\ 
	We obtain
\begin{equation}\label{c.varphi.1eq}
    \begin{array}{lll}
    &\displaystyle\iint\limits_{Q}e^{-2m_0s\alpha}[s^{-1}\xi^{-1}|\Delta\varphi|^2+s\lambda^2\xi|\nabla
    \varphi|^2
    +s^3\lambda^4\xi^3|\varphi|^2]dxdt\\
    &\leq C\Biggl( s^{16}\lambda^{40}\displaystyle\iint\limits_{\omega_0\times(0,T)}
    e^{-8m_0s\hat\alpha+6m_0s\alpha^*}(\hat\xi)^{16}|\varphi|^2dxdt\\
    &\hspace{3cm}
    +s^{15/2}\lambda^{20}\displaystyle\iint\limits_{\mathcal{O}_d\times(0,T)}e^{-4m_0s\hat\alpha+2m_0s
    \alpha^*}(\hat\xi)^{15/2}|\theta|^2dxdt\\
    &\hspace{4cm}+s^{15/2}\lambda^{20}\displaystyle\iint\limits_{Q}e^{-4m_0s\hat\alpha+2m_0s
    \alpha^*}(\hat\xi)^{15/2}|g_1|^2dxdt\Biggr), 
    \end{array}
\end{equation}
    for any $\lambda\geq C$ and $s\geq C(T^5+T^{10})$.\\
    
    Taking into account that $\|\theta\|_{L^2(\Omega)^N}\leq C\|\nabla\times\theta\|_{L^2(\Omega)^{2N-3}}$ and 
    the inequality \eqref{c.intro.inequality.special} with $\varepsilon=\frac{m_0-a_0-1}{m_0+a_0+1}$, 
    $M_1=-\frac{15}{4(m_0+a_0+1)}$ and $M_2=-\frac{10}{(m_0+a_0+1)}$, the second term in the right--hand side 
    of \eqref{c.varphi.1eq} can be 
    estimated by $$\displaystyle\iint\limits_{Q}e^{-2s\alpha^*-2a_0s\alpha^*}|\nabla\times\theta|^2dxdt$$
    and therefore it can be absorbed by the left--hand side of \eqref{c.estimatetheta.end}.\\
    From \eqref{c.estimatetheta.end} and \eqref{c.varphi.1eq} we have
\begin{equation}\label{c.theta.and.varphi}
    \begin{array}{lll}
    &\displaystyle\iint\limits_{Q}e^{-2s\alpha-2a_0s\alpha^*}(s\lambda^2\xi|\nabla(\nabla\times\theta)|^2dxdt
    +s^3\lambda^4\xi^3|\nabla\times\theta|^2)dxdt\\
    &\hspace{1cm}+\displaystyle\iint\limits_{Q}e^{-2m_0s\alpha}[s^{-1}\xi^{-1}|\Delta\varphi|^2
    +s\lambda^2\xi|\nabla\varphi|^2
    +s^3\lambda^4\xi^3|\varphi|^2]dxdt\\
    &\leq Cs^{16}\lambda^{40}\displaystyle\iint\limits_{\omega_0\times(0,T)}
    e^{-8m_0s\hat\alpha+6m_0s\alpha^*}(\hat\xi)^{16}|\varphi|^2dxdt
    +C_{s,\lambda}\|\rho^*(-\ell^{-2}\varphi\chi_{\mathcal{O}}+\gamma^{-2}\varphi)\|^2_{L^2(Q)^N}\\ 
     &\hspace{1.5cm}
     +Cs^{15/2}\lambda^{20}\displaystyle\iint\limits_{Q}e^{-4m_0s\hat\alpha+2m_0s
    \alpha^*}(\hat\xi)^{15/2}|g_1|^2dxdt+C\|\rho^*g_2\|^2_{L^2(Q)^N}\\
    &\hspace{1.5cm}
    +Cs^3\lambda^4\displaystyle\iint\limits_{\omega_0\times (0,T)}e^{-2s\alpha-2a_0s\alpha^*}
     \xi^3|\nabla\times\theta|^2dxdt,
    \end{array}
\end{equation}
    for any $\lambda\geq C$, $s\geq C(T^5+T^{10})$ and $C_{s,\lambda}$ depending on $s,\lambda$.\\
    Choosing $\mathcal{\ell}$ and $\gamma$ large enough 
    ($\mathcal{\ell},\gamma\approx s^{4}\lambda^{5}e^{\lambda\|\eta\|_{\infty}})$, we can absorb the second 
    term in the right--hand side of \eqref{c.theta.and.varphi} by the left--hand side.\\
    Let us estimate the local term concerning $\nabla\times\theta$ in terms of $\varphi$. To do this, we use 
    the
    first equation of \eqref{carleman_adjoint_system} since $\omega\cap\mathcal{O}_d\neq \emptyset$ and
    $\omega_0\subset \mathcal{O}_d$. We have
    $$-(\nabla\times\varphi)_t-\Delta(\nabla\times\varphi)=\nabla\times g_1
    +\mu(\nabla\times\theta),\quad\mbox{in}\,\,\omega_0\times(0,T).$$
    Then,
    \begin{equation*}
    \begin{array}{lll}
    I&:=&s^3\lambda^4\displaystyle\iint\limits_{\omega_0\times(0,T)}e^{-2s\alpha-2a_0s\alpha^*}\xi^3|\nabla\times\theta|^2dxdt\\
    &=& s^3\lambda^4\displaystyle\iint\limits_{\omega_0\times(0,T)}e^{-2s\alpha-2a_0s\alpha^*}\xi^3
    (\nabla\times\theta)(-(\nabla\times\varphi)_t-\Delta(\nabla\times\varphi)-(\nabla\times g_1))dxdt.
    \end{array}
    \end{equation*}
    We introduce an open set $\omega_1\Subset\omega$ such that $\omega_0\Subset\omega_1$ and a positive 
    function $\zeta\in C^2_c(\omega_1)$ such that $\zeta\equiv 1$ in $\omega_0$. Then, after several 
    integration by parts in time and space we have:
\begin{equation*}
    \begin{array}{lll}
    I&=&s^3\lambda^4\displaystyle\iint\limits_{\omega_0\times(0,T)}e^{-2s\alpha-2a_0s\alpha^*}\xi^3
    (\nabla\times\theta)(-(\nabla\times\varphi)_t-\Delta(\nabla\times\varphi)-(\nabla\times g_1))dxdt\\
    &\leq & s^3\lambda^4\displaystyle\iint\limits_{\omega_1\times(0,T)}\zeta\partial_t
    (e^{-2s\alpha-2a_0s\alpha^*}\xi^3)(\nabla\times\theta)(\nabla\times\varphi)dxdt\\
    & &\hspace{1.5cm}
    +s^3\lambda^4\displaystyle\iint\limits_{\omega_1\times(0,T)}e^{-2s\alpha-2a_0s\alpha^*}\xi^3
    ((\nabla\times\theta)_t-\Delta(\nabla\times\theta))(\nabla\times\varphi)dxdt\\
    & &\hspace{1.5cm}-s^3\lambda^4\displaystyle\iint\limits_{\omega_1\times(0,T)}\Delta(\zeta 
    e^{-2s\alpha-2a_0s\alpha^*}\xi^3)
    (\nabla\times\theta)(\nabla\times\varphi)dxdt\\
    & &\hspace{1.5cm}-2s^3\lambda^4\displaystyle\iint\limits_{\omega_1\times(0,T)}\nabla(\zeta 
    e^{-2s\alpha-2a_0s\alpha^*}\xi^3)(\nabla(\nabla\times\theta))(\nabla\times\varphi)dxdt\\
    & &\hspace{1.5cm}-s^3\lambda^4\displaystyle\iint\limits_{\omega_1\times(0,T)}\zeta 
    e^{-2s\alpha-2a_0s\alpha^*}\xi^3(\nabla\times\theta)(\nabla\times g_1)dxdt.
    \end{array}
\end{equation*}
    From the second equation in \eqref{carleman_adjoint_system}, we have that
\begin{equation}\label{aux_I}
    \begin{array}{lll}
    I&\leq  s^3\lambda^4\displaystyle\iint\limits_{\omega_1\times(0,T)}\zeta\partial_t
    (e^{-2s\alpha-2a_0s\alpha^*}\xi^3)(\nabla\times\theta)(\nabla\times\varphi)dxdt\\
    &\hspace{1.5cm} +s^3\lambda^4\gamma^{-2}\displaystyle\iint\limits_{\omega_1\times(0,T)}
    e^{-2s\alpha-2a_0s\alpha^*}\xi^3|\nabla\times\varphi|^2dxdt\\
    &\hspace{1.5cm} -s^3\lambda^4\displaystyle\iint\limits_{\omega_1\times(0,T)}\Delta(\zeta 
    e^{-2s\alpha-2a_0s\alpha^*}\xi^3)
    (\nabla\times\theta)(\nabla\times\varphi)dxdt\\
    &\hspace{1.5cm}-2s^3\lambda^4\displaystyle\iint\limits_{\omega_1\times(0,T)}\nabla(\zeta 
    e^{-2s\alpha-2a_0s\alpha^*}\xi^3)(\nabla(\nabla\times\theta))(\nabla\times\varphi)dxdt\\
    &\hspace{1.5cm}-s^3\lambda^4\displaystyle\iint\limits_{\omega_1\times(0,T)}\zeta 
    e^{-2s\alpha-2a_0s\alpha^*}\xi^3(\nabla\times\theta)(\nabla\times g_1)dxdt\\
    &\hspace{1.5cm} -s^3\lambda^4\displaystyle\iint\limits_{\omega_1\times(0,T)}
    \zeta e^{-2s\alpha-2a_0s\alpha^*}\xi^3(\nabla\times\varphi)(\nabla\times g_2)dxdt.
    \end{array}
\end{equation}    
    Using the estimate
    $$|\partial_t(e^{-2s\alpha-2a_0s\alpha^*}\xi^3)|\leq CTse^{-2s\alpha-2a_0s\alpha^*}(\xi)^{4+1/5},\quad
    \mbox{ for every}\,\, s\geq C$$
    and Young's inequality, we can deduce the following inequalities:
    \begin{equation*}
        \begin{array}{lll}
            I_1&:= s^3\lambda^4\displaystyle\iint\limits_{\omega_1\times(0,T)}\zeta\partial_t
            (e^{-2s\alpha-2a_0s\alpha^*}\xi^3)(\nabla\times\theta)(\nabla\times\varphi)dxdt\\
            &\leq  CTs^4\lambda^4\displaystyle\iint\limits_{\omega_1\times(0,T)}\zeta
            e^{-2s\alpha-2a_0s\alpha^*}\xi^{4+1/5}|\nabla\times\theta||\nabla\times\varphi|dxdt\\
            &\leq  \varepsilon s^3\displaystyle\iint\limits_{\omega_1\times(0,T)}e^{-2s\alpha-2a_0s\alpha^*}
            \xi^3|\nabla\times\theta|^2dxdt\\
            &\quad+C(\varepsilon)s^5\lambda^8\displaystyle\iint\limits_{\omega_1\times(0,T)}
            e^{-2s\alpha-2a_0s\alpha^*}(\xi)^{5+2/5}|\nabla\times\varphi|^2dxdt,
        \end{array}
    \end{equation*}
    for every $s\geq C$ and every $\varepsilon>0$.\\
    Now, using the estimate 
    $$|\Delta(\zeta e^{-2s\alpha-2a_0s\alpha^*}\xi^3 )|\leq Cs^2\lambda^2e^{-2s\alpha-2a_0s\alpha^*}\xi^{5},
    \quad\mbox{ for every}\,\, s\geq C$$
    and again the Young's inequality for the third term in the right--hand side of \eqref{aux_I}, we obtain 
    \begin{equation*}
        \begin{array}{lll}
            I_3&:= -s^3\lambda^4\displaystyle\iint\limits_{\omega_1\times(0,T)}\Delta(\zeta 
            e^{-2s\alpha-2a_0s\alpha^*}\xi^3)(\nabla\times\theta)(\nabla\times\varphi)dxdt\\       
            &\leq  C s^5\lambda^6\displaystyle\iint\limits_{\omega_1\times(0,T)}
            e^{-2s\alpha-2a_0s\alpha^*}\xi^5|\nabla\times\theta||\nabla\times\varphi|dxdt\\
            &\leq  \varepsilon s^3\displaystyle\iint\limits_{\omega_1\times(0,T)}
            e^{-2s\alpha-2a_0s\alpha^*}\xi^3|\nabla\times\theta|^2dxdt
            +C(\varepsilon)s^7\lambda^{12}\displaystyle\iint\limits_{\omega_1\times(0,T)}
            e^{-2s\alpha-2a_0s\alpha^*}\xi^7|\nabla\times\varphi|^2dxdt,
         \end{array}
    \end{equation*}
    for every $s\geq C$ and every $\varepsilon>0$.\\
    Analogously, we can estimate the fourth term in the right--hand side of \eqref{aux_I} by 
    \begin{equation*}
        \begin{array}{lll}
            I_4&:=-2s^3\lambda^4\displaystyle\iint\limits_{\omega_1\times(0,T)}\nabla(\zeta 
            e^{-2s\alpha-2a_0s\alpha^*}\xi^3)(\nabla(\nabla\times\theta))(\nabla\times\varphi)dxdt\\       
            &\leq \varepsilon s\lambda^2\displaystyle\iint\limits_{\omega_1\times(0,T)}
            e^{-2s\alpha-2a_0s\alpha^*}\xi|\nabla(\nabla\times\theta)|^2dxdt
            +C(\varepsilon)s^7\lambda^8\displaystyle\iint\limits_{\omega_1\times(0,T)}
            e^{-2s\alpha-2a_0s\alpha^*}\xi^7|\nabla\times\varphi|^2 dxdt,
         \end{array}
    \end{equation*}
    for every $s\geq C$ and every $\varepsilon>0$.\\
    Additionally, through another integration by part and Young's inequality we can obtain
    \begin{equation*}
        \begin{array}{lll}
            I_5&:=-s^3\lambda^4\displaystyle\iint\limits_{\omega_1\times(0,T)}\zeta 
            e^{-2s\alpha-2a_0s\alpha^*}\xi^3(\nabla\times\theta)(\nabla\times g_1)dxdt\\       
            &\leq \varepsilon\Bigl( s\lambda^2\displaystyle\iint\limits_{\omega_1\times(0,T)}
            e^{-2s\alpha-2a_0s\alpha^*}\xi|\nabla(\nabla\times\theta)|^2dxdt
            +s^3\lambda^4\displaystyle\iint\limits_{\omega_1\times(0,T)}
            e^{-2s\alpha-2a_0s\alpha^*}\xi^3|\nabla\times\theta|^2dxdt\Bigr)\\
            &\hspace{1cm}+C(\varepsilon)s^5\lambda^6\displaystyle\iint\limits_{Q}
            e^{-2s\alpha-2a_0s\alpha^*}\xi^5|g_1|^2dxdt,
         \end{array}
    \end{equation*}
    for every $s\geq C$ and every $\varepsilon>0$.\\
    \begin{equation*}
        \begin{array}{lll}
    	I_6&:=-s^3\lambda^4\displaystyle\iint\limits_{\omega_1\times(0,T)}
    	\zeta e^{-2s\alpha-2a_0s\alpha^*}\xi^3(\nabla\times\varphi)(\nabla\times g_2)dxdt\\
    	&\leq C \Bigl(\|\rho^*g_2\|_{L^2(Q)^N}
    	+s^7\lambda^{12}\displaystyle\iint\limits_{\omega_1\times (0,T)}e^{-2s\alpha-2a_0s\alpha^*}
    	\xi^7|\nabla\varphi|^2 dxdt\\
    	&\hspace{1cm} +s^6\lambda^8\displaystyle\iint\limits_{\omega_1\times (0,T)}e^{-4s\alpha-2a_0s\alpha^*}
    	 \xi^6|\nabla\times(\nabla\times\varphi)|^2 dxdt \Bigr).
       \end{array}
    \end{equation*}
    
    We use Lemma \ref{lemma.integral.weights} in the Appendix in order to obtain an appropriate upper bound for the 
    last term in the   right--hand side on the previous inequality 
    $$s^7\lambda^{12}\displaystyle\iint\limits_{\omega_1\times (0,T)}e^{-2s\alpha-2a_0s\alpha^*}
    \xi^7|\nabla\varphi|^2 dxdt\quad\mbox{ and }\quad 
    \varepsilon s^{-1}\displaystyle\iint\limits_{Q}e^{-2m_0s\alpha}\xi^{-1}|\Delta\varphi|^2dxdt,$$
    for every $\varepsilon>0$. \\
    
    Putting together \eqref{c.theta.and.varphi} and the previous estimates, we have 
\begin{equation}\label{c.theta.and.varphi2}
    \begin{array}{lll}
    &\displaystyle\iint\limits_{Q}e^{-2s\alpha-2a_0s\alpha^*}(s\lambda^2\xi|\nabla(\nabla\times\theta)|^2dxdt
    +s^3\lambda^4\xi^3|\nabla\times\theta|^2)dxdt\\
    &\hspace{1cm}+\displaystyle\iint\limits_{Q}e^{-2m_0s\alpha}[s^{-1}\xi^{-1}|\Delta\varphi|^2
    +s\lambda^2\xi|\nabla\varphi|^2
    +s^3\lambda^4\xi^3|\varphi|^2]dxdt\\
    &\leq Cs^{16}\lambda^{40}\displaystyle\iint\limits_{\omega_0\times(0,T)}
    e^{-8m_0s\hat\alpha+6m_0s\alpha^*}(\hat\xi)^{16}|\varphi|^2dxdt
    +C\displaystyle\iint\limits_{Q}e^{-2a_0s\alpha^*}|g_2|^2dxdt\\ 
     &\hspace{1.5cm}
     +Cs^{15/2}\lambda^{20}\displaystyle\iint\limits_{Q}e^{-4m_0s\hat\alpha+2m_0s
    \alpha^*}(\hat\xi)^{15/2}|g_1|^2 dxdt\\
    &\hspace{1.5cm}
    +Cs^7\lambda^{12}\displaystyle\iint\limits_{\omega_1\times (0,T)}e^{-2s\alpha-2a_0s\alpha^*}
     \xi^7|\nabla\times\varphi|^2 dxdt,
    \end{array}
\end{equation}
    for any $\lambda\geq C$, $s\geq CT^{10}$.\\
    
    On the other hand, considering  open sets $\omega_2,\omega_3\Subset\omega$ such that 
    $\omega_1\Subset \omega_2\Subset \omega_3\subset \omega$, we can deduce that 
\begin{equation}
    \begin{array}{lll}
     &s^7\lambda^{12}\displaystyle\iint\limits_{\omega_1\times (0,T)}e^{-2s\alpha-2a_0s\alpha^*}
     \xi^7|\nabla\times\varphi|^2dxdt\\
     &\leq C(\varepsilon)s^{15}\lambda^{24}\displaystyle\iint\limits_{\omega_3\times(0,T)}
     e^{2(m_0-2)s\alpha^*-4a_0s\alpha^*}(\hat\xi)^{15}|\varphi|^2dxdt\\
     &\hspace{1cm}
     +\varepsilon\Biggl(\displaystyle\iint\limits_{Q}e^{-2m_0s\alpha}[s^{-1}\xi^{-1}|\Delta\varphi|^2
        +s\lambda^2\xi|\nabla\varphi|^2+s^3\lambda^4\xi^3|\varphi|^2]dxdt\Biggr),
    \end{array}
\end{equation}
    for any $\lambda\geq C$, $s\geq CT^{10}$ and any $\varepsilon>0$.\\
    Taking into account that $m_0>a_0+1$, there 
    exists a constant $C>0$ such that
    \begin{equation}\label{c.end.proof}
    s^{16}\lambda^{40}\displaystyle\iint\limits_{\omega_0\times(0,T)}e^{-8sm_0\hat\alpha+6sm_0\alpha^*}
    (\hat\xi)^{16}|\varphi|^2dxdt\leq Cs^{15}\lambda^{24}\displaystyle\iint\limits_{\omega_3\times(0,T)}
    e^{2(m_0-2)s\alpha^*-4a_0s\alpha^*}(\hat\xi)^{15}|\varphi|^2dxdt.
    \end{equation}
    From \eqref{c.theta.and.varphi2}--\eqref{c.end.proof}, we conclude the proof of Proposition 
    \ref{carleman_carleman1}.\\            
\subsection{\normalsize Null controllability of the linear system}
    In this section we are concerned in the null controllability of the linear coupled Stokes system
\begin{equation}\label{3.2.main.system}
    \left\{
    \begin{aligned}
        \begin{array}{llll}
        y_{t}-\Delta y+\nabla p=f_1+h1_{\omega}+(-\ell^{-2}\chi_{\mathcal{O}}+\gamma^{-2})q &\text{ in }& Q,\\
        -z_{t}-\Delta z+\nabla \pi=f_2+\mu (y-y_d)\chi_{\mathcal{O}_d} &\text{ in }& Q,\\
        \nabla \cdot y=0, \nabla \cdot z=0  &\text{ in }& Q, \\
        y=z=0&\text{ on }&\Sigma, \\
        y(\cdot,0)=y_0(\cdot),\quad z(\cdot,T)=0& \text{ in }&\Omega,
        \end{array}
    \end{aligned}
    \right.
\end{equation}
    where the functions $f_1$ and $f_2$ are in appropriate weighted spaces. We look for a control 
    $h\in L^2(\omega\times(0,T))^N$ such that, under 
    suitable properties on $f_1, f_2$, the solution to \eqref{3.2.main.system} satisfies 
    $y(\cdot,T)=0$ in $\Omega$.\\
    To do this, let us first state a Carleman inequality with weight functions not vanishing in $t=0$. Let 
    $\tilde{\ell}\in C^1([0,T])$ be a positive function in $[0,T)$ such that: 
    \begin{equation*}
    \tilde{\ell}(t)=T^2/4\ \forall 
  	t\in [0,T/2]\  \text{ and }\ \tilde{\ell}(t)=t(T-t)\ \forall  t\in [T/2,T].
  	\end{equation*} Now, we introduce the following weight
    functions
\begin{equation}\label{3.2.carleman.weights}
    \begin{aligned}
    &\beta(x,t) = \dfrac{e^{12\lambda\|\eta\|_{\infty}}-e^{\lambda(10\|\eta\|_{\infty}+\eta(x))}} 
    {\tilde{\ell}^{5}(t)},\quad \tau(x,t)=\dfrac{e^{\lambda(10\|\eta\|_{\infty}+\eta(x))}}{\tilde{\ell}^{5}
    (t)}
    ,\\
    &\beta^*(t) = \max_{x\in\overline{\Omega}} \beta(x,t),\quad \quad
     \tau^*(t) = \min_{x\in\overline{\Omega}} \tau(x,t),\\
    &\widehat\beta(t) = \min_{x\in\overline{\Omega}} \beta(x,t),\quad 
    \quad \,\,\,\,\widehat\tau(t) = \max_{x\in\overline{\Omega}} \tau(x,t).
    \end{aligned}
\end{equation}

\begin{lema}\label{3.2.lemma.carleman2}
    Let $s$ and $\lambda$ like in Proposition \ref{carleman_carleman1}. Then, there exists a constant $C>0$
    (depending on $s,\lambda,\omega,\mathcal{O},T,\ell,\gamma$ and $\mu$) such that every solution 
    $(\varphi,\theta)$ of \eqref{carleman_adjoint_system} satisfies
    \begin{equation}\label{3.2.ine.carleman2}
    \begin{array}{lll}
        &\|\varphi(\cdot,0)\|^2_{L^2(Q)^N} 
        +\displaystyle\iint\limits_{Q}e^{-2m_0s\beta^*}(\tau^*)^{3}|\varphi|^2dxdt
        +\displaystyle\iint\limits_{Q}e^{-2(a_0+1)s\beta^*}(\tau^*)^3|\theta|^2dxdt\\
        &\leq C\Biggl(\displaystyle\iint\limits_{Q}e^{-2a_0s\beta^*}(\hat\tau)^{15}|g_1|^2dxdt
        +\displaystyle\iint\limits_{Q}e^{-2a_0s\beta^*}|g_2|^2dxdt\\
        &\hspace{3cm}+\displaystyle\iint\limits_{\omega\times(0,T)}e^{-4a_0s\beta^*+2(m_0-2)s\beta}
        (\hat\tau)^{15}|\varphi|^2dxdt\Biggr).
    \end{array}
    \end{equation} 
\end{lema}
	
    \noindent\textit{Proof of Lemma \ref{3.2.lemma.carleman2}.} By construction 
    $\alpha=\beta$ and $\xi=\tau$
    in $\Omega\times(T/2,T)$,  so that
    \begin{equation*}
    \begin{array}{ll}
    	&\displaystyle\int\limits_{T/2}^T\displaystyle\int\limits_{\Omega}(e^{-2(a_0+1)s\alpha^*}(\xi^*)^3|\theta|^2 +
    e^{-2sm_0\alpha^*}(\xi^*)^3|\varphi|^2)dxdt\\
    &=\displaystyle\int\limits_{T/2}^T\displaystyle\int\limits_{\Omega}(e^{-2(a_0+1)s\beta^*}(\tau^*)^3|\theta|^2
    +e^{-2sm_0\beta^*}(\tau^*)^3|\varphi|^2)dxdt.	
    \end{array}
	\end{equation*}

    Therefore, it follows from Proposition \ref{carleman_carleman1} the estimate
    \begin{equation*}
    \begin{array}{lll}
        &\displaystyle\int\limits_{T/2}^T\int\limits_{\Omega}(e^{-2(a_0+1)s\beta^*}(\tau^*)^3|\theta|^2
        +e^{-2sm_0\beta^*}(\tau^*)^3|\varphi|^2)dxdt\\
        &\leq C\Biggl( \displaystyle\iint\limits_{Q}e^{-2a_0s\alpha^*}(\hat\xi)^{5}|g_1|^2dxdt
        +\displaystyle\iint\limits_{Q}e^{-2a_0s\alpha^*}|g_2|^2dxdt\\
        &\hspace{3cm}+\displaystyle\iint\limits_{\omega\times(0,T)}e^{-4a_0s\alpha^*+2(m_0-2)s\alpha}
        (\hat\xi)^{15}|\varphi|^2dxdt\Biggr). 
    \end{array}
    \end{equation*}     
    Since $\tilde\ell(t)=t(T-t)$ for any $t\in [T/2,T]$ and 
    $$e^{-2a_0s\beta^*}\geq C,\quad e^{-2a_0s\beta^*}(\tau^*)^{5}\geq C\quad\mbox{and}\quad
    e^{-4a_0s\beta^*+2(m_0-2)s\beta}(\hat\tau)^{15}\geq C\,\, \mbox{in}\,\, [0,T/2],$$
    we readily get           
    \begin{equation}\label{carleman2.estimate.right}
    \begin{array}{lll}
        &\displaystyle\int\limits_{T/2}^T\int\limits_{\Omega}(e^{-2(a_0+1)s\beta^*}(\tau^*)^3|\theta|^2
        +e^{-2sm_0\beta^*}(\tau^*)^3|\varphi|^2)dxdt\\
        &\leq C\Biggl( \displaystyle\iint\limits_{Q}e^{-2a_0s\beta^*}(\hat\tau)^{5}|g_1|^2dxdt
        +\displaystyle\iint\limits_{Q}e^{-2a_0s\beta^*}|g_2|^2dxdt\\
        &\hspace{4cm}+\displaystyle\iint\limits_{\omega\times(0,T)}e^{-4a_0s\beta^*+2(m_0-2)s\beta}
        (\hat\tau)^{15}|\varphi|^2dxdt\Biggr). 
    \end{array}
    \end{equation}   
    
    Now, we introduce a function $\nu\in C^1([0,T])$ such that $\nu\equiv 1$ in $[0, T/2],\, \nu\equiv 0$ in
    $[3T/4,T]$. It is easy to see that $(\nu\varphi,\nu\pi_{\varphi})$ and $(\nu\theta,\nu\pi_{\varphi})$
    satisfies the system
\begin{equation}\label{carleman2.system.nu.thetavarphi}
    \left\{
    \begin{aligned}
        \begin{array}{llll}
        -(\nu\varphi)_{t}-\Delta (\nu\varphi)+\nabla(\nu\pi_{\varphi}) 
        =\nu(g_1+\mu\theta\chi_{\mathcal{O}_d})-\nu'\varphi &\text{ in }& Q,\\
        (\nu\theta)_{t}-\Delta (\nu\theta)+\nabla (\nu\pi_{\theta})
        =\nu(g_2-\ell^{-2}\varphi\chi_{\mathcal{O}}+\gamma^{-2}\varphi)+\nu'\theta &\text{ in }& Q,\\
        \nabla \cdot (\nu\varphi)=\nabla \cdot (\nu\theta)=0, \nabla \cdot q=0  &\text{ in }& Q, \\
        \nu\varphi=\nu\theta=0&\text{ on }&\Sigma, \\
        (\nu\varphi)(T)=0,\quad \nu\theta(0)=0& \text{ in }&\Omega,
        \end{array}
    \end{aligned}
    \right.
\end{equation}     
    Using classical energy estimate for both $\nu\varphi$ and $\nu\theta$, which solve the Stokes system 
    \eqref{carleman2.system.nu.thetavarphi} we get
    \begin{equation*}
    \begin{array}{lll}
    &\|\varphi(0)\|^2_{L^2(Q)^N}+\|\varphi\|^2_{L^2(0,T/2;H_0^1(\Omega)^N)}\\
    &\leq C\Biggl( \displaystyle\frac{1}{T^2}\|\varphi\|^2_{L^2(T/2,T/4;L^2(\Omega)^N)}
    +\|\theta\|^2_{L^2(0,3T/4;L^2(\mathcal{O}_d)^N)} +\|g_1\|^2_{L^2(0,3T/2;L^2(\Omega)^N)}\Biggr)
    \end{array}
    \end{equation*}
    and 
    \begin{equation*}
    \begin{array}{lll}
    \|\theta\|^2_{L^2(0,T/2;H_0^1(\Omega)^N)}&
    \leq C\Biggl( \displaystyle\frac{1}{T^2}\|\theta\|^2_{L^2(T/2,3T/4;L^2(\Omega)^N)}\\
    &\hspace{1cm}+\|\nu(-\ell^{-2}\varphi\chi_{\mathcal{O}}
    +\gamma^{-2}\varphi)\|^2_{L^2(0,3T/4;L^2(\Omega)^N)}
    +\|g_2\|^2_{L^2(0,3T/2;L^2(\Omega)^N)}\Biggr).
    \end{array}
    \end{equation*}
    Taking into account that 
    $$e^{-2sm_0\beta^*}(\tau^*)^3\geq C>0 \quad e^{-2(a_0+1)s\beta^*}(\tau^*)^3\geq C>0,\quad \forall 
    t\in [T/2,3T/4]$$
    and 
    $$e^{-2a_0s\beta^*}(\hat\tau)^{5}\geq C>0\quad 
    e^{-2a_0s\beta^*}>e^{-4a_0s\beta^*}\geq C>0,\quad\forall t\in [0,3T/4],$$
    we have 
    \begin{equation}\label{carleman2.estimate.left1}
    \begin{array}{lll}
    &\|\varphi(0)\|^2_{L^2(\Omega)^N}+\displaystyle\int\limits_{0}^{T/2}\int\limits_{\Omega}
    e^{-2m_0s\beta^*}(\tau^*)^{3}|\varphi|^2dxdt
    +\displaystyle\int\limits_{0}^{T/2}\int\limits_{\Omega}e^{-2(a_0+1)s\beta^*}(\tau^*)^{3}|\theta|^2dxdt\\
    &\leq C\Biggl(\displaystyle\int\limits_{T/2}^{3T/2}\int\limits_{\Omega}
    [e^{-2m_0s\beta^*}(\tau^*)^{3}|\varphi|^2
    +e^{-2(a_0+1)s\beta^*}(\tau^*)^{3}|\theta|^2]dxdt\\
    &\hspace{1cm}
    +\|\nu\mu e^{-2a_0s\beta^*}\theta\|^2_{L^2(0,3T/4;L^2(\mathcal{O}_d)^N)}
    +\|\nu(-\ell^{-2}\varphi\chi_{\mathcal{O}}
    +\gamma^{-2}\varphi)\|^2_{L^2(0,3T/4;L^2(\Omega)^N)}\\
    &\hspace{1cm}+\displaystyle\int\limits_{0}^{3T/4}\int\limits_{\Omega}
    \Bigl[e^{-2a_0s\beta^*}(\tau^*)^{5}|g_1|^2+e^{-2a_0s\beta^*}|g_2|^2\Bigr]dxdt
    \Biggr). 
    \end{array}
    \end{equation}
	Thus, from \eqref{carleman2.estimate.right} and \eqref{carleman2.estimate.left1} we have at this moment
	\begin{equation}\label{carleman2.estimate.right.left}
    \begin{array}{lll}
        &\|\varphi(0)\|^2_{L^2(\Omega)^N}+\displaystyle\iint\limits_{Q}(e^{-2(a_0+1)s\beta^*}(\tau^*)^3|\theta|^2
        +e^{-2sm_0\beta^*}(\tau^*)^3|\varphi|^2)dxdt\\
        &\leq C\Biggl( \displaystyle\iint\limits_{Q}e^{-2a_0s\beta^*}(\hat\tau)^{5}|g_1|^2dxdt
        +\displaystyle\iint\limits_{Q}e^{-2a_0s\beta^*}|g_2|^2dxdt\\
        &\hspace{0.5cm}+\displaystyle\iint\limits_{\omega\times(0,T)}e^{-4a_0s\beta^*+2(m_0-2)s\beta}
        (\hat\tau)^{15}|\varphi|^2dxdt\\
        &\hspace{0.5cm}+\|\nu\mu e^{-2a_0s\beta^*}\theta\|^2_{L^2(0,3T/4;L^2(\mathcal{O}_d)^N)}
    	+\|\nu(-\ell^{-2}\varphi\chi_{\mathcal{O}}+\gamma^{-2}\varphi)\|^2_{L^2(0,3T/4;L^2(\Omega)^N)}
        \Biggr). 
    \end{array}
    \end{equation}     
    
    Observe that if $\ell$ and $\gamma$ are large enough, the last term in the right--hand side of
    \eqref{carleman2.estimate.right.left} can be absorbed by the left--hand side. In addition, considering 
    $\theta^*(x,t)=e^{-2s\beta^*}\theta(x,t)$ instead $\nu\theta$ in \eqref{carleman2.system.nu.thetavarphi} 
    and using standard energy estimate for the system associated to $\theta$, we obtain
    \begin{equation}\label{carleman2.estimate.thetalocal.right}
    \begin{array}{lll}
        \displaystyle\int\limits_{0}^{3T/4}\int\limits_{\Omega}\nu^2\mu^2e^{-4a_0s\beta^*}|\theta|^2dxdt
        &\leq C\Biggl( \displaystyle\iint\limits_{Q}e^{-4a_0s\beta^*}|g_2|^2dxdt
        +\displaystyle\frac{1}{\ell^4}\int\limits_{0}^T\int\limits_{\mathcal{O}_d}e^{-4a_0s\beta^*}
        |\varphi|^2dxdt\\
        &\hspace{0.2cm}+\displaystyle\frac{1}{\gamma^4}\iint\limits_{Q}e^{-4a_0s\beta^*}|\varphi|^2dxdt
        +\displaystyle\iint\limits_{Q}e^{-4a_0s\beta^*}(\tau^*)^{6/5}|\theta|^2dxdt
        \Biggr). 
    \end{array}
    \end{equation}  
    Putting together \eqref{carleman2.estimate.right.left}, \eqref{carleman2.estimate.thetalocal.right} and
    taking again $\ell$ and $\gamma$ large enough, we obtain the desired inequality  
    \eqref{3.2.ine.carleman2}.\\
    
    \begin{Obs}\label{obs.carleman2.add.term}
        In order to establish a null controllability result for the system \eqref{3.2.main.system} with 
        suitable weight functions, see Theorem \ref{prop.null.control}, observe that on the left--hand side
        of \eqref{carleman_obs_ineq1} it is possible to add the term
         $$\displaystyle\iint\limits_{Q}e^{-4a_0s\beta^*}(\tau^*)^{3}|\theta|^2dxdt.$$
         This is a consequence of the inequalities $a_0+1<m_0\leq 2a_0,\,\, a_0\geq \frac{5}{4}$. 
    \end{Obs}

    Now, we are ready to prove the null controllability of system \eqref{3.2.main.system}. The idea is to
    look a solution in an appropriate weighted functional space. Let us introduce the following space\\
    \begin{equation}
    \begin{array}{lll}
    	E&:=\{(y,z,\pi_y,\pi_z,h):e^{a_0s\beta^*}(\hat\tau)^{-5/2}y\in L^2(Q)^N, 
    	e^{a_0s\beta^*}z\in L^2(Q)^N,\\
    	&\hspace{1cm}e^{2a_0s\beta^*-(m_0-2)s\hat\beta}(\hat\tau)^{-15/2}h1_{\omega}\in L^2(Q)^N,\\
  		&\hspace{1cm}e^{a_0s\beta^*}(\hat\tau)^{-15/2}y\in L^2(0,T;H^2(\Omega)^N)\cap L^{\infty}(0,T;V),\\
    	&\hspace{1cm}e^{a_0s\beta^*}(\tau^*)^{-c_0}z\in L^2(0,T;H^2(\Omega)^N)\cap L^{\infty}(0,T;V),\,\,
    	c_0\geq\frac{5}{2},\\
    	&\hspace{1cm}
    	e^{m_0s\beta^*}(\hat\tau^*)^{-3/2}(y_t-\Delta y+\nabla\pi_y-(-\ell^{-2}\chi_{\mathcal{O}}+\gamma^{-2})z
    	-h1_{\omega})\in L^2(Q)^N,\\
    	&\hspace{1cm}
    	e^{2a_0s\beta^*}(\hat\tau^*)^{-3/2}(-z_t-\Delta z+\nabla\pi_z
    	-\mu(y-y_d)\chi_{\mathcal{O}_d})\in L^2(Q)^N
    	\}.
    \end{array}
	\end{equation}      
	It is clear that $E$ is a Banach space for the following norm:
	\begin{equation*}
    \begin{array}{lll}
    	\|(y,z,\pi_y,\pi_z,h)\|_E:=&\|e^{a_0s\beta^*}(\hat\tau)^{-5/2}y\|_{L^2(Q)^N}+ 
    	\|e^{a_0s\beta^*}z\|_{L^2(Q)^N}\\
    	&+\|e^{2a_0s\beta^*-(m_0-2)s\hat\beta}(\hat\tau)^{-15/2}h1_{\omega}\|_{L^2(Q)^N}\\
  		&+\|e^{a_0s\beta^*}(\hat\tau)^{-15/2}y\|_{L^2(0,T;H^2(\Omega)^N)}
  		+\|e^{a_0s\beta^*}(\hat\tau)^{-15/2}y\|_{L^{\infty}(0,T;V)}\\
    	&+\|e^{a_0s\beta^*}(\tau^*)^{-c_0}z\|_{L^2(0,T;H^2(\Omega)^N)}
    	+\|e^{a_0s\beta^*}(\tau^*)^{-c_0}z\|_{L^{\infty}(0,T;V)}\\
    	&+
    	\|e^{m_0s\beta^*}(\tau^*)^{-3/2}(y_t-\Delta y+\nabla\pi_y-(-\ell^{-2}\chi_{\mathcal{O}}+\gamma^{-2})z
    	-h1_{\omega})\|_{L^2(Q)^N}\\
    	&+
    	\|e^{2a_0s\beta^*}(\tau^*)^{-3/2}(-z_t-\Delta z+\nabla\pi_z
    	-\mu(y-y_d)\chi_{\mathcal{O}_d})\|_{L^2(Q)^N}.
    \end{array}
	\end{equation*} 
\begin{Obs}
	Observe in particular that $(y,z,\pi_y,\pi_z,h)\in E$ implies $y(\cdot, T)=0$ in $\Omega$.
\end{Obs}
	
\begin{teo}\label{prop.null.control}
	Assume the hypothesis of Lemma \ref{3.2.lemma.carleman2} and 
	\begin{equation}\label{condition.f1_f2}
	y_0\in V,\,\,\,\, e^{m_0s\beta^*}(\tau^*)^{-3/2}f_1\in L^2(Q)^N,\,\,\,\,
	e^{2a_0s\beta^*}(\tau^*)^{-3/2}f_2\in L^2(Q)^N.
	\end{equation}
	Then, we can find a control $h\in L^2(\omega\times(0,T))^N$ such that the associated solution 
	$(y,z,\pi_y,\pi_z,h)$ to \eqref{3.2.main.system} satisfies $(y,z,\pi_y,\pi_z,h)\in E$.
\end{teo}
\noindent\textit{Proof of Theorem \ref{prop.null.control}.} 
    Let us introduce the following constrained extremal problem:
	 \begin{equation}\label{prop.null.control.extremal.problem}
	 \left\{
	 \begin{array}{lll}
	 \inf & \left\{
	\begin{array}{lll}
		&\displaystyle\frac{1}{2}\Bigl(\iint\limits_{Q}e^{2a_0s\beta^*}(\hat\tau)^{-5}|y|^2dxdt
	 +\iint\limits_{Q}e^{2a_0s\beta^*}|z|^2dxdt\\
	 &\hspace{1cm}
	 +\displaystyle\iint\limits_{\omega\times(0,T)}e^{4a_0s\beta^*-2(m_0-2)s\hat\beta}
	 (\hat\tau)^{-15}|h|^2dxdt\Bigr)
	\end{array}\right.\\
	 &\mbox{subject to}\,\, h\in L^2(Q),\,\, supp\, h\subset \omega\times (0,T),\,\,\mbox{and}\,\,\\
     &\begin{aligned}\left\{
        \begin{array}{llll}
        y_{t}-\Delta y+\nabla\pi_y=f_1+h\chi_{\omega}+(-\ell^{-2}\chi_{\mathcal{O}}+\gamma^{-2})z 
        &\text{ in }& Q,\\
        -z_{t}-\Delta z+\nabla\pi_z=f_2+\mu (y-y_d)\chi_{\mathcal{O}_d} &\text{ in }& Q,\\
        \nabla \cdot y=0, \nabla \cdot z=0  &\text{ in }& Q, \\
        y=z=0&\text{ on }&\Sigma, \\
        y(\cdot,0)=y_0(\cdot),\quad y(\cdot,T)=0,\quad z(\cdot,T)=0& \text{ in }&\Omega.
        \end{array}
        \right.
    &\end{aligned}
    &\end{array}\right.
	 \end{equation}
	Assume that this problem admits a unique solution $(\hat y, \hat z, \hat{\pi_y},\hat{\pi_z},\hat h)$. Then, 
	from the Lagrange's principle there exists dual variables 
	$(\hat\varphi, \hat\theta, \hat{\pi}_{\varphi},\hat{\pi}_{\theta})$ such that
	\begin{equation}\label{prop.nullcontrol.dualvariables}
	\begin{array}{llll}
	&\hat y= e^{-2a_0s\beta^*}(\hat\tau)^{5}
	(-\hat\varphi_t-\Delta\hat\varphi +\nabla\hat\pi_{\varphi}-\mu\hat\theta_{\mathcal{O}_d})&\mbox{in}& Q,\\
	&\hat z=e^{-2a_0s\beta^*}(\hat\theta_t-\Delta\hat\theta+\nabla\hat\pi_{\theta}
	-(-\ell^{-2}\chi_{\mathcal{O}}+\gamma^{-2})\hat\varphi) &\mbox{in}& Q, \\
	&\hat h=e^{-4a_0s\beta^*+2(m_0-2)s\hat\beta}(\hat\tau)^{15}\hat\varphi &\mbox{in}& Q,\\
	&\hat y=\hat z=0 &\mbox{on}& \Sigma.
	\end{array}
	\end{equation}
	Now, following the arguments established in \cite{FCGIP04}, we introduce the space $P_0$ of functions 
	$(y,z,\pi_y,\pi_z)\in C^2(\overline{Q})^{2N+2}$ such that
	\begin{enumerate}
	\item [i)] $\nabla\cdot y=\nabla\cdot z=0$\,\, in\, $Q$.
	\item [ii)] $y=z=0$\,\, on\,\, $\Sigma$.
	\item [iii)] $\displaystyle\int\limits_{\omega_0}\pi_{\varphi} dx=0$.
	\end{enumerate}
	We also consider the bilinear form $a(\cdot,\cdot)$ over $P_0\times P_0$ defined by:
	\begin{multline*} 
	a((\hat\varphi,\hat\theta,\hat\pi_{\varphi},\hat\pi_{\theta}),(w,z,\pi_w,\pi_z))=:\\
	   \displaystyle\iint\limits_{Q}e^{-2a_0s\beta^*}(\hat\tau)^{5}
	    (-\hat\varphi_t-\Delta\hat\varphi +\nabla\hat\pi_{\varphi}-\mu\hat\theta_{\mathcal{O}_d})
	    (-y_t-\Delta y +\nabla\pi_{y}-\mu z_{\mathcal{O}_d})\, dxdt\\
	 +\displaystyle\iint\limits_{Q}e^{-2a_0s\beta^*}(\hat\theta_t-\Delta\hat\theta+\nabla\hat\pi_{\theta}
	    -(-\ell^{-2}\chi_{\mathcal{O}}+\gamma^{-2})\hat\varphi)
	    (z_t-\Delta z+\nabla\pi_{z})\\
	  -\displaystyle\iint\limits_{Q}e^{-2a_0s\beta^*}(\hat\theta_t-\Delta\hat\theta
	    +\nabla\hat\pi_{\theta}-(-\ell^{-2}\chi_{\mathcal{O}}+\gamma^{-2})\hat\varphi)
	    (\ell^{-2}\chi_{\mathcal{O}}+\gamma^{-2}) w)\, dxdt\\
	 + \displaystyle\iint\limits_{\omega\times (0,T)}\!\!\!e^{-4a_0s\beta^*+2(m_0-2)s\hat\beta}
	 (\hat\tau)^{15}
	    \hat\varphi w\, dxdt   ,
	\end{multline*}

    for every $(w,z,\pi_w,\pi_z)\in P_0$,  and a linear form
    \begin{equation}\label{prop.null.control.def.lineal.form}
        \langle G, (w,z,\pi_w,\pi_z)\rangle:=\displaystyle\iint\limits_{Q}f_1\cdot w\, dxdt+
        \displaystyle\iint\limits_{Q}f_2\cdot z\, dxdt
        +\displaystyle\int\limits_{\Omega}y_0(\cdot)\cdot w(\cdot,0)\, dx. 
    \end{equation}
    Taking into account this definitions, one can see that, if the functions $\hat y,\hat z$ and $\hat h$ solve
    \eqref{prop.null.control.extremal.problem}, we must have 
    \begin{equation}\label{prop.null.control.identity.bilineal.lineal}
        a((\hat\varphi,\hat\theta,\hat\pi_{\varphi},\hat\pi_{\theta}),(w,z,\pi_w,\pi_z))
        =\langle G, (w,z,\pi_w,\pi_z)\rangle, \quad \forall  (w,z,\pi_w,\pi_z)\in P_0. 
    \end{equation}
    Observe that Carleman inequality \eqref{3.2.ine.carleman2} holds for all $(w,z,\pi_w,\pi_z)\in P_0$. 
    Consequently, 
    \begin{equation}\label{prop.nullcontrol.estimate1.bilinealform}
        \begin{array}{lll}
        &\displaystyle\iint\limits_{Q}e^{-2m_0s\beta^*}(\tau^*)^{3}|z|^2 dxdt
        +\displaystyle\iint\limits_{Q}e^{-2(a_0+1)s\beta^*}(\tau^*)^3|w|^2 dxdt\\
        &+\displaystyle\iint\limits_{Q}e^{-2a_0s\beta^*}(\tau^*)^3|w|^2 dxdt
        +\|w(0)\|^2_{L^2(\Omega)^N}\leq C a((w,z,\pi_w,\pi_z),(w,z,\pi_w,\pi_z)),
        \end{array}     
    \end{equation}
    for every $(w,z,\pi_w,\pi_z)\in P_0$.\\
    Therefore, $a(\cdot,\cdot):P_0\times P_0\longmapsto\mathbb{R}$  is symmetric, definite positive bilinear 
    form
    on $P_0$. We denote by $P$ the completion of $P_0$ for the norm induced by $a(\cdot,\cdot)$. Then,
    $a(\cdot,\cdot)$ is well--defined, continuous and again definite positive on $P$. Furthermore, in view
    of the Carleman inequality \eqref{3.2.ine.carleman2}, the assumption \eqref{condition.f1_f2} and 
    \eqref{prop.nullcontrol.estimate1.bilinealform}, the linear form 
    $(w,z,\pi_w,\pi_z)\longmapsto \langle G, (w,z,\pi_w,\pi_z)\rangle$  is well--defined and continuous on
    $P$. Indeed, for every $(w,z,\pi_w,\pi_z)\in P$, 
    \begin{equation*}
    \begin{array}{lll}
        \langle G, (w,z,\pi_w,\pi_z)\rangle  
        &\leq  \|e^{(a_0+1)s\beta^*}(\tau^*)^{-3/2}f_1\|_{L^2(Q)^N}
        \|e^{-(a_0+1)s\beta^*}(\tau^*)^{3/2}w\|_{L^2(Q)^N}\\
        &\hspace{0.5cm}+
        \|e^{m_0s\beta^*}(\tau^*)^{-3/2}f_2\|_{L^2(Q)^N}\|e^{-m_0s\beta^*}(\tau^*)^{3/2}z\|_{L^2(Q)^N}
        +\|y_0\|_{H}\|w(0)\|_{H}\\
        &\leq  \|e^{m_0s\beta^*}(\tau^*)^{-3/2}f_1\|_{L^2(Q)^N}
        \|e^{-(a_0+1)s\beta^*}(\tau^*)^{3/2}w\|_{L^2(Q)^N}\\
        & \hspace{0.5cm}
        +\|e^{2a_0s\beta^*}(\tau^*)^{-3/2}f_2\|_{L^2(Q)^N}\|e^{-m_0s\beta^*}(\tau^*)^{3/2}z\|_{L^2(Q)^N}
        +\|y_0\|_{H}\|w(0)\|_{H}.
    \end{array}
    \end{equation*} 
    Using \eqref{prop.nullcontrol.estimate1.bilinealform} and the density of $P_0$ in $P$, we find
    \begin{equation*}
        \begin{array}{ll}
        &\langle G, (w,z,\pi_w,\pi_z)\rangle
        \leq C\Bigl( \|e^{m_0s\beta^*}(\tau^*)^{-3/2}f_1\|_{L^2(Q)^N}\\
        &\hspace{3cm}
        +\|e^{2a_0s\beta^*}(\tau^*)^{-3/2}f_2\|_{L^2(Q)^N}+\|y_0\|_{H}\Bigr)\|(w,z,\pi_w,\pi_z)\|
        _{P}.
        \end{array}    
    \end{equation*}
    Hence, from Lax--Milgram's Lemma, there exists a unique 
    $(\hat\varphi,\hat\theta,\hat\pi_{\varphi},\hat\pi_{\theta})\in P$ satisfying:
    \begin{equation}\label{prop.nullcontrol.estimate2.identity}
        a((\hat\varphi,\hat\theta,\hat\pi_{\varphi},\hat\pi_{\theta}),(w,z,\pi_w,\pi_z))
        =\langle G, (w,z,\pi_w,\pi_z)\rangle, \quad \forall  (w,z,\pi_w,\pi_z)\in P.     
    \end{equation}
    Let us set $(\hat y, \hat z,\hat h)$ like in \eqref{prop.nullcontrol.dualvariables} and remark that 
    $(\hat y,\hat z,\hat\pi_{y},\hat\pi_z,\hat h)$ verifies
    \begin{equation*}
    \begin{array}{lll}
      a((\hat\varphi,\hat\theta,\hat\pi_{\varphi},\hat\pi_{\theta}),
      (\hat\varphi,\hat\theta,\hat\pi_{\varphi},\hat\pi_{\theta}))=& 
      \displaystyle\iint\limits_{Q}e^{2a_0s\beta^*}(\hat\tau)^{-5}|\hat y|^2 dxdt
      +\displaystyle\iint\limits_{Q}e^{2a_0s\beta^*}|\hat z|^2 dxdt\\
      & +\displaystyle\iint\limits_{\omega\times(0,T)}e^{4a_0s\beta^*-2(m_0-2)s\hat\beta}(\hat\tau)^{-15}
      |\hat h| dxdt< +\infty.  
    \end{array}
    \end{equation*}
    
    Let us prove that $(\hat y,\hat z)$ is, together with some $(\hat\pi_y,\hat\pi_z)$, the weak solution 
    of the Stokes system in \eqref{prop.null.control.extremal.problem} for $h=\hat h$. In fact, we introduce 
    the (weak) solution $(\tilde y,\tilde z,\tilde\pi_y,\tilde\pi_z)$ to the Stokes system
    \begin{equation}\label{prop.nullcontrol.system.aux1}
    \left\{
        \begin{array}{llll}
        \tilde y_{t}-\Delta\tilde y+\nabla\tilde\pi_y=f_1+\hat h1_{\omega}
        +(-\ell^{-2}\chi_{\mathcal{O}}+\gamma^{-2})\tilde z &\text{ in }& Q,\\
        -\tilde z_{t}-\Delta\tilde z+\nabla\tilde\pi_z=f_2+\mu(\tilde y-\tilde y_d)\chi_{\mathcal{O}_d} 
        &\text{ in }& Q,\\
        \nabla \cdot \tilde y=0, \nabla \cdot \tilde z=0  &\text{ in }& Q, \\
        \tilde y=\tilde z=0&\text{ on }&\Sigma, \\
        \tilde y(\cdot,0)=y_0(\cdot),\quad \tilde z(\cdot,T)=0& \text{ in }&\Omega,
        \end{array}
    \right.
    \end{equation}
    Clearly, $(\tilde y,\tilde z)$ is the unique solution of \eqref{prop.nullcontrol.system.aux1} defined
    by transposition. This means that, for every $(a,b)\in L^2(Q)^{2N}$, 
    \begin{equation}\label{prop.nullcontrol.transposition1}
        \langle (\tilde y,\tilde z), (a,b)\rangle_{L^2(Q)^N}=\langle y_0,\varphi(0)\rangle_{L^2(\Omega)}
        +\langle (f_1+\hat h1_{\omega},f_2), (\varphi,\theta)\rangle_{L^2(Q)^N},
    \end{equation}
    where $(\varphi,\theta)$ is, together with some $(\pi_\varphi,\pi_\theta)$, the solution to
    \begin{equation}\label{prop.nullcontrol.system.aux2}
    \left\{
        \begin{array}{llll}
        L^*(\varphi,\theta)=(a,b) &\text{ in }& Q,\\
        \nabla\cdot \varphi=0, \nabla \cdot \theta=0  &\text{ in }& Q, \\
        \varphi=\theta=0&\text{ on }&\Sigma, \\
        \varphi(\cdot,T)=0,\quad \theta(\cdot,0)=0& \text{ in }&\Omega,
        \end{array}
    \right.
    \end{equation}
    and $L^*$ is the adjoint operator of $L$  given by:
    $$L(\tilde y,\tilde z):=(\tilde y_{t}-\Delta\tilde y+\nabla\tilde\pi_y
    -(-\ell^{-2}\chi_{\mathcal{O}}+\gamma^{-2})\tilde z,
    -\tilde z_{t}-\Delta\tilde z+\nabla\tilde\pi_z-\mu(\tilde y-\tilde y_d)\chi_{\mathcal{O}_d}).$$
    From \eqref{prop.nullcontrol.dualvariables} and \eqref{prop.null.control.identity.bilineal.lineal},
    we see that $(\hat y,\hat z)$ also satisfies \eqref{prop.nullcontrol.transposition1}. Consequently, 
    $(\hat y,\hat z)=(\tilde y,\tilde z)$ and $(\hat y,\hat z)$ is, together with some 
    $(\hat\pi_y,\hat\pi_z)=(\tilde\pi_y,\tilde\pi_z)$, the weak solution to the system 
    \eqref{prop.nullcontrol.system.aux1}.\\
    Finally, we must see that $(\hat y,\hat z,\hat\pi_y,\hat\pi_z,\hat h)\in E$. We already know that 
    $$e^{a_0s\beta^*}(\hat\tau)^{-5/2}\hat y,\,\, e^{a_0s\beta^*}\hat z,\,\,\,
    e^{2a_0s\beta^*-(m_0-2)s\hat\beta}(\hat\tau)^{-15/2}\hat h 1_{\omega}\in L^2(Q)^N$$
    and (see hypothesis \eqref{condition.f1_f2})
    $$e^{m_0s\beta^*}(\tau^*)^{-3/2}f_1\in L^2(Q)^N\quad\mbox{and}\quad 
	e^{2a_0s\beta^*}(\tau^*)^{-3/2}f_2\in L^2(Q)^N.$$
	Thus, it only remains to check that
	$$e^{a_0s\beta^*}(\hat\tau)^{-15/2}\hat y,\,\, e^{a_0s\beta^*}(\tau^*)^{-c_0} z
	\in L^2(0,T;H^2(\Omega)^N)\cap L^{\infty}(0,T;V),$$
	where $c_0\geq \frac{5}{2}.$
	
	\begin{enumerate}
	\item We define the functions
	    \begin{equation*}
	    \begin{array}{llll}
	    y^*:=e^{a_0s\beta^*}(\hat\tau)^{-15/2}\hat y, & & z^*:=e^{a_0s\beta^*}(\tau^*)^{-c_0}\hat z\\
	    \pi^*_{y}:=e^{a_0s\beta^*}(\hat\tau)^{-15/2}\hat\pi_y, & & \pi^*_z:=e^{a_0s\beta^*}(\tau^*)^{-c_0}
	    \hat\pi_z
	    \end{array}
	    \end{equation*}
        and
        \begin{equation*}
        \begin{array}{llll}
        f^*_1:=e^{a_0s\beta^*}(\hat\tau)^{-15/2}(f_1+h1_{\omega}), & &
        z^{**}:=e^{a_0s\beta^*}(\hat\tau)^{-15/2}(-\ell^{-2}\chi_{\mathcal{O}}+\gamma^{-2})z\\
        f^*_2:=e^{a_0s\beta^*}(\tau^*)^{-c_0} f_2, & & 
        y^{**}:=e^{a_0s\beta^*}(\tau^*)^{-c_0}(y-y_d)\chi_{\mathcal{O}_d}. 
        \end{array}
        \end{equation*}
        Then $(y^*,\pi^*_y,z^*,\pi^*_z)$ satisfies:
        \begin{equation}\label{prop.nullcontrol.system.aux3}
        \left\{
        \begin{array}{llll}
        y^*_{t}-\Delta y^*+\nabla\pi^*_y=f^*_1+z^{**}+(e^{3/2s\beta^*}(\hat\tau)^{-15/2})^{'}\hat y 
        &\text{ in }& Q,\\
        -z^*_{t}-\Delta z^*+\nabla\pi^*_z=f^*_2+y^{**}+(e^{1/2s\beta^*}(\hat\tau)^7)^{'}\hat z 
        &\text{ in }& Q,\\
        \nabla \cdot y^*=0, \nabla \cdot z^*=0  &\text{ in }& Q, \\
        y^*=z^*=0&\text{ on }&\Sigma, \\
        y^*(\cdot,0)=e^{3/2s\beta^*(0)}(\hat\tau(0))^{-15/2}y_0(\cdot),
        \quad z^*(\cdot,T)=0& \text{ in }&\Omega,
        \end{array}
    	\right.
    	\end{equation}
    
    \item Now, we prove that the right--hand side of the main equations in \eqref{prop.nullcontrol.system.aux3}
        is in $L^2(Q)^N.$
        \begin{itemize}
        \item $|e^{a_0s\beta^*}(\hat\tau)^{-15/2}f_1|\leq C e^{a_0s\beta^*}|\hat\tau|^{-15/2}|f_1|
                \leq Ce^{m_0s\beta}|\tau^*|^{-3/2}|f_1|.$
        \item $|e^{a_0s\beta^*}(\hat\tau)^{-15/2}h1_{\omega}|
                \leq Ce^{2a_0s\beta^*-(m_0-2)s\beta^*}(\hat\tau)^{-15/2}|h|1_{\omega}.$
        \item $|z^{**}|=|e^{a_0s\beta^*}(\hat\tau)^{-15/2}(-\ell^{-2}\chi_{\mathcal{O}}+\gamma^{-2})z|
        		\leq Ce^{a_0s\beta^*}|\hat{z}|.$
		\item $|(e^{3/2s\beta^*}(\hat\tau)^{-15/2})^{'}\hat y|
				\leq Cse^{a_0s\beta^*}|\tau^*|^{6/5}|\hat{y}|
				\leq C e^{a_0s\beta^*}|\hat{\tau}|^{-5/2}|\hat{y}|$.
        \item $|f_2^{**}|=|e^{a_0s\beta^*}(\tau^*)^{-c_0} f_2|\leq Ce^{(a_0+1)s\beta^*}|\tau^*|^{-c_0}|f_2|$.
        \item $|(e^{1/2s\beta^*}(\hat\tau)^7)^{'}\hat z|\leq Ce^{a_0s\beta^*}|\hat{z}|$.
        \item $|y^{**}|=|e^{a_0s\beta^*}(\tau^*)^{-c_0}(y-y_d)\chi_{\mathcal{O}_d}|
        		\leq Ce^{a_0s\beta^*}|\hat{\tau}|^{-5/2}|\hat{y}|+Ce^{a_0s\beta^*}|\tau^*|^{-c_0}|y_d|.$
        \end{itemize}
        Observe that $y^{**}\in L^2(Q)^N$ thanks to the hypothesis \eqref{hypothesis.weight.mainresult}.
      \end{enumerate}
    Taking into account $a)-b)$ and  $y_0\in V$, we have $y^*,z^*\in L^2(0,T;H^2(\Omega)^N)\cap L^{\infty}(0,T;V)$
    (see Lemma \ref{carleman_regularityH3} in Section 3.1).\\
    This concludes the proof of Proposition \ref{prop.null.control}.
	\begin{Obs}
	Before starting the last section, it is important to consider small data in order to prove our main result,
	Theorem \ref{th_mainresult}. Thus, we impose that
	\begin{equation}\label{delta_condition}
	\|f_1\|_{L^2(Q)^N}+\|f_2\|_{L^2(Q)^N}+ \|y_0\|_{V}\leq \delta,
	\end{equation}
	where $\delta$ is a small positive number.
	\end{Obs}

\section{\normalsize Proof of the main result} In this section we give the proof of Theorem \ref{th_mainresult}
	throughout classical arguments such like in \cite{FCGIP04}. The results obtained in the previous section allow
	us to locally invert a nonlinear operator associated to the nonlinear system
 		\begin{equation*}
 		\left\{
        \begin{array}{llll}
        y_{t}-\Delta y+(y\cdot\nabla)y+\nabla\pi_y=h1_{\omega}+(\ell^{-2}\tilde{\chi}_{\mathcal{O}}+\gamma^{-2})z 
        &\text{ in }& Q,\\
        -z_{t}-\Delta z+(z,\nabla^t)y-(y,\nabla)z+\nabla\pi_z=\mu(y-y_d)\chi_{\mathcal{O}_d} 
        &\text{ in }& Q,\\
        \nabla \cdot y=0, \nabla \cdot z=0  &\text{ in }& Q, \\
        y=z=0&\text{ on }&\Sigma, \\
        y(\cdot,0)=y_0(\cdot),\quad z(\cdot,T)=0& \text{ in }&\Omega.
        \end{array}
    	\right.
    	\end{equation*}
    to do this, we will apply an inverse function theorem of the Luisternik's kind
    \cite{hamilton1982inverse}, which will allow us to complete the proof of theorem \ref{th_mainresult}.
     More precisely, we will use the following theorem.
    \begin{teo}\label{teo_inverse_mapping}
	Suppose that $\mathcal{B}_1,\mathcal{B}_2$ are Banach spaces and 
	$$\mathcal{A}:\mathcal{B}_1\to \mathcal{B}_2$$
	is a continuously differentiable map. We assume that for $b_1^0\in \mathcal{B}_1, b_2^0\in \mathcal{B}_2$ the
	equality
	\begin{equation}\label{b_1^0}
	\mathcal{A}(b_1^0)=b_2^0
	\end{equation}
	holds and $\mathcal{A}'(b_1^0):\mathcal{B}_1\to \mathcal{B}_2$ is an epimorphism.Then there exists 
	$\delta >0$ such that for any $b_2\in \mathcal{B}_2$ which satisfies the condition 
	$$\|b_2^0-b_2\|_{\mathcal{B}_2}<\delta$$
	there exists a solution $b_1\in \mathcal{B}_1$ of the equation$$\mathcal{A}(b_1)=b_2.$$
	\end{teo}
	
	We apply this theorem for the spaces $\mathcal{B}_1:=E$ and 
	$$\mathcal{B}_2:=\{(f_1,f_2,y_0)\in X_1^*\times X_2^*\times V: f_1,f_2,y_0\,\, \mbox{satisfies}
	\,\,\,\eqref{delta_condition}\},$$
	where $X_1^*:= L^2(e^{m_0s\beta^*}(\tau^*)^{-3/2}(0,T);L^2(\Omega)^N)$ and 
	$X_2^*:=L^2(e^{2a_0s\beta^*}(\tau^*)^{-3/2}(0,T);L^2(\Omega)^N)$.\\
	We define the operator $\mathcal{A}$ by the formula
	\begin{equation}
	\begin{array}{lll}
		\mathcal{A}(y,z,\pi_y,\pi_z,h)&:= (y_{t}-\Delta y+(y\cdot\nabla)y+\nabla\pi_y
		-(\ell^{-2}\tilde{\chi}_{\mathcal{O}}+\gamma^{-2})z-h1_{\omega},\\
		 & \hspace{0.8cm}
		 -z_{t}-\Delta z+(z,\nabla^t)y-(y,\nabla)z+\nabla\pi_z-\mu(y-y_d)\chi_{\mathcal{O}_d},y(\cdot,0)),
	\end{array}
	\end{equation}
	for every $(y,z,\pi_y,\pi_z,h)\in \mathcal{B}_1$.\\
	
	Let us see that $\mathcal{A}$ is of class $C^1(\mathcal{B}_1, \mathcal{B}_2)$. Indeed, notice that all the
	terms in $\mathcal{A}$ are linear, except for $(y\cdot\nabla)y$ and $z,\nabla^t)y-(y,\nabla)z$, then, we only
	have to check that these nonlinear terms are well--defined and depend continuously on the data. Thus, we will
	prove that the bilinear operator 
	$((y^1,z^1,\pi^1_y,\pi^1_z,h^1),(y^2,z^2,\pi^2_y,\pi^2_z,h^2))\longmapsto (y^1\cdot\nabla)y^2$	
	is continuous from $\mathcal{B}_1\times\mathcal{B}_1$ to $X_1^*$, and the bilinear forms 
	$((y^1,z^1,\pi^1_y,\pi^1_z,h^1),(y^2,z^2,\pi^2_y,\pi^2_z,h^2))\longmapsto (y^1\cdot\nabla)z^2$,\,\,
	$((y^1,z^1,\pi^1_y,\pi^1_z,h^1),(y^2,z^2,\pi^2_y,\pi^2_z,h^2))\longmapsto (z^1\cdot\nabla^t)y^2$
	are continuous from $\mathcal{B}_1\times\mathcal{B}_1$ to $X_2^*$.\\
	In fact, notice that (see the definition of the space $E$):
	$$e^{a_0s\beta^*}(\hat\tau)^{-15/2}y\in L^2(0,T;L^{\infty}(\Omega)^N)$$ and
	$$\nabla (e^{a_0s\beta^*}(\hat\tau)^{-15/2}y)\in L^{\infty}(0,T;L^2(\Omega)^{N\times N}).$$
	Consequently, we obtain
	\begin{equation*}
	\begin{aligned}
		&\|e^{m_0s\beta^*}(\tau^*)^{-3/2}(y^1\cdot\nabla)y^2\|_{L^2(Q)^N}\\
		&\leq C \|(e^{a_0s\beta^*}(\hat\tau)^{-15/2}y^1\cdot\nabla) e^{a_0s\beta^*}
		(\hat\tau)^{-15/2}y^2\|_{L^2(Q)^N}\\
		&\leq C\|e^{2s\beta^*}(\hat\tau)^{-15/2}y^1\|_{L^2(0,T;L^\infty(\Omega)^N)}
		\|e^{a_0s\beta^*}(\hat\tau)^{-15/2}y^2\|_{L^{\infty}(0,T;V)}.
	\end{aligned}
	\end{equation*}
	On the other hand, for $c_0\geq 5/2$, 
	$$e^{a_0s\beta^*}(\tau^*)^{-c_0}z\in L^2(0,T;L^{\infty}(\Omega)^N)$$ and
	$$\nabla (e^{a_0s\beta^*}(\tau^*)^{-c_0}z)\in L^{\infty}(0,T;L^2(\Omega)^{N\times N}).$$
	Then,
	\begin{equation*}
	\begin{aligned}
		&\|e^{2a_0s\beta^*}(\tau^*)^{-3/2}(y^1\cdot\nabla)z^2\|_{L^2(Q)^N}\\
		&\leq C\|e^{a_0s\beta^*}(\hat\tau)^{-15/2}y^1\|_{L^2(0,T;L^\infty(\Omega)^N)}
		\|e^{a_0s\beta^*}(\tau^*)^{-c_0}z^2\|_{L^{\infty}(0,T;V)},
	\end{aligned}
	\end{equation*}
	and analogously,
	\begin{equation*}
	\begin{aligned}
		&\|e^{2a_0s\beta^*}(\tau^*)^{-3/2}(z^1\cdot\nabla)y^2\|_{L^2(Q)^N}\\
		&\leq C\|e^{a_0s\beta^*}(\tau^*)^{-c_0}z^1\|_{L^2(0,T;L^\infty(\Omega)^N)}
		\|e^{a_0s\beta^*}(\hat\tau)^{-15/2}y^2\|_{L^{\infty}(0,T;V)}.
	\end{aligned}
	\end{equation*}
	Notice that $\mathcal{A}'(0,0,0):\mathcal{B}_1\to \mathcal{B}_2$ is given by
	$$(y_{t}-\Delta y+\nabla\pi_y-(\ell^{-2}\tilde{\chi}_{\mathcal{O}}+\gamma^{-2})z-h1_{\omega},
	-z_{t}-\Delta z+\nabla\pi_z-\mu(y-y_d)\chi_{\mathcal{O}_d}, y(\cdot,0)),$$
	for all\,\,$ (y,z,\pi_y,\pi_z,h)\in \mathcal{B}_1.$\\
	In virtue of Theorem \ref{prop.null.control}, this functional satisfies 
	$Im (\mathcal{A}'(0,0,0))=\mathcal{B}_2$.\\
	Let $b_1^0=(0,0,0)$ and $b_2^0=(0,0)$. Then equation \eqref{b_1^0} obviously holds.
	So all necessary conditions to apply Theorem \ref{teo_inverse_mapping} are 
	fulfilled. Therefore there exists a positive number $\delta$ such that, if 
	$\|y(\cdot,0)\|_{V}\leq\delta$, we can find a control $h\in L^2(\omega\times(0,T))^N$ 
	and an associated solution $(y,z,\pi_y,\pi_z)$ to \eqref{1.1.main_system} satisfying $y(\cdot,T)=0$ in
 	$\Omega$. This finishes the proof of Theorem \ref{th_mainresult}.
 	
\section{\normalsize Conclusion and open problems}
	In this article, we mentioned the main results on robust control for the 
	$N$--dimensional Navier--Stokes system with Dirichlet boundary conditions. These results has 
	also allowed us to characterise the follower control $v$ and its disturbance function $\psi$ through
	a nonlinear coupled system. Once this step has finished, we used the  robust pair $(v,\psi)$ to prove 
	the null controllability of the leader control $h$. The main novelties are the Carleman inequalities
	for coupled Stokes system, which involves new relationships between the weight functions and the
	robustness parameters $\ell, \gamma$, see Proposition 
	\eqref{carleman_carleman1} and Lemma \ref{relationship.weights}. To conclude, we present now some open
	problems arising from our study: 
	\begin{itemize}
	\item If instead of considering in  the hierarchical strategy a zero objective for the leader control 
		$h$ in \eqref{1.1.main_system}, the objective may be a trajectory $(\overline{y},\overline{\pi})$ of the
		uncontrolled system:
		\begin{equation*}
 		\left\{
        \begin{array}{llll}
        \overline{y}_{t}-\Delta \overline{y}+(\overline{y}\cdot\nabla)\overline{y}+\nabla\overline{\pi}=0       
        &\text{ in }& Q,\\
        \nabla \cdot \overline{y}=0, &\text{ in }& Q, \\
        \overline{y}=0&\text{ on }&\Sigma, \\
        \overline{y}(\cdot,0)=\overline{y}_0(\cdot) & \text{ in }&\Omega,
        \end{array}
    	\right.
    	\end{equation*}
    	 So we may ask if is it possible to prove the local exact controllability to  trajectories of system 
    	 \eqref{1.1.main_system}. That is, does there exist a control $h$ such that for the corresponding
    	  solution to \eqref{robust_coupled_system_linear}  satisfies $y(T)=\overline y (T)$?
    \item Some  null controllability results for the $N$--dimensional Navier--Stokes equation
     (\cite{carreno2013local}, \cite{coron2014local}) allow to act on the system by means of few controls. 
    	Is it possible to
    	extend these results to a robust Stackelberg strategy? Is is possible to ask 
   	  the leader control $h$ to have one vanishing component?	 
	\item  Is it possible to extend the  results in this paper to  Navier--slip boundary conditions?. In 
		other words, can we say something about the existence and uniqueness of saddle points for the 
		Navier--Stokes system with Navier--slip conditions? Does we have the null controllability for
		the leader control $h$?  
	\item Finally, it would be interesting to study the problems proposed in this paper
		 to other models such as water waves (Korteweg--de Vries equation), interaction fluid--heat 
		 (Boussinesq system), micropolar fluids, models of turbulence, among others.
	\end{itemize}
	    
\section{\normalsize Appendix: some technical results}
	In the following results, it will be assumed that $N=2$ or $N=3$.\\
	From the relation between $\alpha^*$ and $\hat\alpha$, it is possible to prove the following 
	inequality.
\begin{lema}\label{relationship.weights}
	For any $\varepsilon>0$, any $M_1,M_2\in \mathbb{R}$, there exists $\lambda_0>0$ and
  	$C=C(\varepsilon,M_1, M_2)>0$ such that 
	\begin{equation}\label{c.intro.inequality.special}
    	e^{s\alpha^*}\leq Cs^{M_1}\lambda^{M_2}(\hat\xi)^{M_1}e^{s(1+\varepsilon)\hat\alpha}
	\end{equation}
    for every $\lambda>\lambda_0$.\\ 
\end{lema}
	\noindent\textit{Proof of Lemma \ref{relationship.weights}}. Recall that 
	$$\alpha^*(t) := \max_{x\in\overline{\Omega}} \alpha(x,t),\quad \quad
	\widehat\alpha(t) := \min_{x\in\overline{\Omega}} \alpha(x,t)\quad \mbox{and}\quad
    \widehat\xi(t) := \max_{x\in\overline{\Omega}} \xi(x,t).$$
    From the definition of $\alpha^*$ and $\widehat{\alpha}$,
     $\widehat{\alpha}(t)=F(\lambda)\alpha^*$, where 
    $F(\lambda):=\displaystyle\frac{e^{2\lambda\|\eta_0\|_{\infty}}-e^{\lambda\|\eta_0\|_{\infty}}}
    {e^{2\lambda\|\eta_0\|_\infty}-1}$. It is easy to check that $F(\lambda)\to 1$ to 
    $\lambda\to+\infty$ and
    $F(\lambda)\to 1/2$ to $\lambda\to 0^+$. Additionally, by construction of $F(\lambda)$, for any 
    $\varepsilon>0$, there exists $\lambda_0>0$ such that, for every $\lambda\geq \lambda_0$
    $$F(\lambda)+\varepsilon F(\lambda)>1.$$
    In consequence, exists a positive constant $C=C(\varepsilon,M_1,\tilde{M_2})$ such that the inequality 
    $$\lambda^{\tilde {M_2}}e^{(1-(1+\varepsilon)F(\lambda))s\alpha^*}\leq Cs^{M_1}(\hat{\xi})^{M_1}$$
    holds for any $M_1,\tilde{M_2}\in\mathbb{R}$.\\
    This completes the proof of Lemma \ref{relationship.weights}.\\
     
	As a consequence of Lemma \ref{relationship.weights}, for $a_0\geq 2$ and $m_0$ satisfying  
	$a_0< m_0\leq a_0+2 $, we can deduce the next result:
\begin{lema}\label{lemma.integral.weights}
	Under the hypothesis of Lemma \ref{relationship.weights}, for any $\omega\Subset\Omega$ and any  $u\in V$,
	there exists 
	$\lambda_0>0$ and
  	$C=C(\varepsilon,\tilde M_1, \tilde M_2)>0$ such that 
  	$$s^{\tilde M_1}\lambda^{\tilde M_2}\displaystyle\iint\limits_{\omega\times(0,T)}
  	e^{-4s\hat\alpha-2a_0s\alpha^*}(\hat{\xi})^{\tilde M_1}|\Delta u|^2dxdt
  	\leq Cs^{-1}\displaystyle\iint\limits_{Q}e^{-2m_0 s\alpha^*}(\hat{\xi})^{-1}|\Delta u|^2dxdt.$$ 
\end{lema}
	\noindent\textit{Sketch of the proof}. Taking 
	$\varepsilon=\displaystyle\frac{2}{m_0-a_0}-1, \, M_1=-\displaystyle\frac{\tilde{M_1}+1}{2(m_0-a_0)}$ and 
	$M_2=-\displaystyle\frac{\tilde{M_2}+1}{2(m_0-a_0)}$ in \eqref{c.intro.inequality.special}, the proof is
	direct.\\

\end{document}